\documentclass[red,11pt,a4paper]{article}

\usepackage{cite}

\usepackage{amsmath}
\usepackage{amscd}
\usepackage{amssymb}
\usepackage{latexsym}
\usepackage{xcolor}
\usepackage{eufrak}
\numberwithin{equation}{section}
\usepackage[normalem]{ulem}

\def\CA{{\mathcal A}}

\def\CD{{\mathcal D}}
\def\CS{{\mathcal S}}

\newcommand{\CL}{{\mathcal L}}
\newcommand{\CC}{{\mathcal C}}
\newcommand{\CR}{{\mathcal R}}
\newcommand{\CF}{{\mathcal F}}
\newcommand{\CB}{{\mathcal B}}
\newcommand{\CT}{{\mathcal T}}

\newcommand{\CX}{{\mathcal X}}
\newcommand{\C}{\mathbb{C}}
\newcommand{\kv}{\vk_\vv}
\newcommand{\kw}{\vk_\vw}
\newcommand{\kbv}{\vk_{\bar \vv}}
\newcommand{\ds}{\,{\rm d}s}
\newcommand{\ebt}{e^{\beta t}}

\newcommand{\fr}{\mathfrak{r}}
\newcommand{\fq}{\mathfrak{q}}
\newcommand{\fz}{\mathfrak{z}}

\newcommand{\lpqt}{L_p(0,T;L_q(\Omega))}
\newcommand{\linf}{L_\infty(\Omega\times(0,T))}
\newcommand{\linfsp}{L_\infty(\Omega)}

\renewcommand{\div}{{\rm div}\,}

\def\d{\partial}

\newcommand{\de}{\partial}

\newcommand{\vO}{\mathbf{O}}
\newcommand{\bvr}{\bar \vr}

\newcommand{\bzeta}{\bar \fz}

\renewcommand{\d}{\partial}
\newcommand{\bv}{\boldsymbol{v}}

\newcommand{\N}{\boldsymbol{n}}

\newcommand{\bI}{\boldsymbol{I}}

\newcommand{\pd}{\partial}

\newcommand{\bV}{{\bold V}}
\newcommand{\vk}{{\bold k}}

\newcommand{\lap}{\Delta}

\newcommand{\Div}{\operatorname{div}}

\newcommand{\vf}{{\vc f}}
\newcommand{\vg}{{\vc g}}

\newcommand{\vS}{{\bf S}}
\newcommand{\vw}{{\bf w}}

\newcommand{\vr}{\varrho}

\newcommand{\vu}{\vc{u}}
\newcommand{\vv}{\vc{v}}

\newcommand{\vc}[1]{{\bf #1}}

\newcommand{\qed}{\rightline{ $\square$}}
\newcommand{\Grad}{\nabla}

\newcommand{\pt}{\partial_{t}}
\newcommand{\ptb}[1]{\partial_{t}(#1)}

\newcommand{\dt}{{\rm d} t }

\newcommand{\lr}[1]{\left( #1 \right)}

\newcommand{\eq}[1]{\begin{equation}
\begin{split}
#1
\end{split}
\end{equation}}
\newcommand{\eqh}[1]{\begin{equation*}
\begin{split}
#1
\end{split}
\end{equation*}}
\newcommand{\ep}{\varepsilon}

\newcommand{\R}{\mathbb{R}}
\newtheorem{thm}{Theorem}[section]
\newtheorem{lem}[thm]{Lemma}
\newtheorem{prop}[thm]{Proposition}
\newtheorem{df}{Definition}
\newtheorem{rmk}{Remark}
\newtheorem{cor}[thm]{Corollary}

\title{Maximal Regularity  for  Compressible Two-Fluid System }

\author{Tomasz Piasecki\thanks{Institute of Applied Mathematics and Mechanics, University of Warsaw, Banacha 2, 02-097 Warsaw, Poland, E-mail: \texttt{t.piasecki@mimuw.edu.pl}},\quad Ewelina Zatorska\thanks{Imperial College London, Department of Mathematics,  London SW7 2AZ, United Kingdom. Email: \texttt{e.zatorska@imperial.ac.uk}}}
\date{}

\begin{document}
 \maketitle 

\centerline{{\em Dedicated to Professor Yoshihiro Shibata on the occasion of his 70th birthday.}}

\vskip10mm

\abstract
We investigate a compressible two-fluid Navier-Stokes type system with a single velocity field and algebraic closure for the pressure law. The constitutive relation involves densities of both fluids through an implicit function. We are interested in regular solutions in a $L_p-L_q$ maximal regularity setting. We show that such solutions exists locally in time and, under additional smallness assumptions on the initial data, also globally. Our proof rely on appropriate transformation of the original problem, application of Lagrangian coordinates and maximal regularity estimates for associated linear problem.


\medskip


\section{Introduction}

\smallskip 

    In the present paper, we analyze a bi-fluid compressible system. It arises, in particular, as a result of averaging procedure applied to liquid-gas free boundary problem (see for instance \cite{Ishii}).  
    We assume a common velocity field and pressure for both fluids (the so-called {\it algebraic pressure closure}), for justification of this reduction and presentation of related models we refer to \cite{Bresch_chapter}. \\
    Our system of equations reads:
\begin{subequations}\label{spec}
\begin{align}
&\ptb{\alpha^\pm\vr^\pm}+\Div(\alpha^{\pm}\vr^\pm\vu)=0,\label{spec_1}\\
&\ptb{(\alpha^{+}\vr^++\alpha^-\vr^-)\vu}+\Div((\alpha^{+}\vr^++\alpha^-\vr^-)\vu\otimes\vu)-\Div \vS(\vu)+\Grad p=\vc{0},\\
&\alpha^++\alpha^-=1,\\
& p = p^+ = p^- \label{spec_4}.
\end{align}
\end{subequations}

By $p^+$, $p^-$ we denote the internal barotropic pressures for each fluid with the explicit forms:
\eq{\label{CP}
p^+=\big(\vr^+\big)^{\gamma_+}, \qquad p^-=\big(\vr^-\big)^{\gamma-}, 
}
where $\gamma^\pm>1$ are given constants, and we assume without loss of generality that
\begin{equation} \label{a:gamma}
\gamma^+ \geq \gamma^-.
\end{equation}

The stress tensor obeys the {\it Newton rheological law}
	\begin{equation}\label{chF:Stokes}
	\vS(\vu)= 2\mu\mathbf{D}(\vu)+\nu\Div \vu{\bf{I}},
	\end{equation}
where $\mathbf{D}(\vu)=\frac{1}{2}\left(\Grad \vu+(\Grad \vu)^{T}\right)$ and $\mu$ and $\nu$ are the nonnegative viscosity coefficients. 

We consider the system \eqref{spec} on a domain $\Omega\subset \R^3$ supplied with the Dirichlet boundary condition for the velocity
\eq{\label{av}
\vu|_{\d\Omega}=0 \quad \mbox{ for } \  t\in (0,T),}
and with the initial conditions
\eq{\label{Ic_1}
\alpha^+ \vr^+ \vert_{t=0} = R_0, \quad \alpha^- \vr^- \vert_{t=0} = Q_0, \qquad 
\vu|_{t=0}=\vu_0.
}
We assume that the total initial density is separated away from zero, i.e.
\eq{\label{init:pos} R_0 \ge 0, \quad Q_0 \ge 0, \quad R_0+Q_0\geq\kappa \;\; \textrm{for some} \; \kappa>0} 
and that the initial data 
$$\alpha^\pm \vert_{t=0} = \alpha^\pm_0,\qquad \vr^\pm \vert_{t=0} = \vr^\pm_0$$
satisfy the following compatibility conditions
\eq{\label{Ic_2}
\alpha^+_0 + \alpha^-_0 = 1, \qquad \alpha^\pm_0 \ge 0, \qquad
p^+(\vr^+_0) = p^-(\vr^-_0).
}

We immediately switch to a reformulation of the system \eqref{spec}.  Introducing the notation
$$R =\vr^+\alpha^+,\qquad Q =\vr^-\alpha^-, \qquad Z= \vr^+,\qquad \alpha=\alpha^+,$$
we check that the pressure $p$ is expressed in terms of $R,Q$. 
For this purpose we observe that, by \eqref{spec_4} and \eqref{CP},
$$
(\vr^-\alpha^-)^{\gamma^-}=(\alpha^-)^{\gamma^-}(\vr^+)^{\gamma^+}
= \left[\frac{\vr^+(1-\alpha)}{\vr^+}\right]^{\gamma^-}(\vr^+)^{\gamma^+}.
$$
Therefore we have
\begin{equation}\label{pZ}
 p=P(R,Q)= {Z}^{\gamma_+},
\end{equation}
for $Z=Z(R,Q)$ such that 
\eq{\label{TZ}
Q = \lr{1-\frac{R}{Z} } Z^\gamma,\quad  \mbox{with}\quad \gamma=\frac{\gamma_+}{\gamma_-},
}
and
\eq{\label{RleqZ}R\leq Z.}
Note that by \eqref{RleqZ} we immediately deduce that $0\leq \alpha\leq 1$.\\ 

The system \eqref{spec}-\eqref{av} can be therefore transformed to the following form
\begin{subequations}\label{S}
\begin{align}
  &\pt R+ \Div(R\vu )=0,\label{SR}\\[3pt] 
  &\pt Q + \Div(Q\vu )=0 ,\label{ST}\\[3pt]
&\ptb{(R+Q)\vu}+\Div((R+Q)\vu\otimes\vu)-\Div \vS(\vu) +\Grad Z^{\gamma^+}=\vc{0},\label{SM}\\[3pt]
&Z=Z(R,Q),\\[3pt]
&R|_{t=0}=R_0, \quad Q|_{t=0}=Q_0, \quad \vu|_{t=0}=\vu_0,\\[3pt]
&\vu|_{\d\Omega}=0,
\end{align}
\end{subequations} 
where $Z$ is related to $R$ and $Q$ through the non-explicit formula \eqref{TZ}. 

Note that by \eqref{spec_4},\eqref{CP} and \eqref{pZ} the triplet $(R,Q,Z(R,Q))$ uniquely determines $(\vr^+,\vr^-,\alpha^+,\alpha^-)$. Therefore it is legitimate to consider the system \eqref{S} instead of \eqref{spec}.

\subsection{Notation}
We use standard notation $L_p(\Omega)$ and $W^k_p(\Omega)$ for Lebesgue and Sobolev spaces, respectively. Furthermore by $L_p(I;X)$, where $I \subset \R_+$ and $X$ is a Banach space, we denote a Bochner space.  

Next, we recall that for $0<s<\infty$ and $m$ a smallest integer larger than $s$ we define Besov spaces as intermediate spaces
\begin{equation} \label{def:bsqp0} 
B^{s}_{q,p}(\Omega)=(L_q(\Omega),W^m_q(\Omega))_{s/m,p},
\end{equation}
where $(\cdot,\cdot)_{s/m,p}$ is the real interpolation functor, see \cite[Chapter 7]{Ad}. In particular,
\begin{equation} \label{def:bsqp}
B^{2(1-1/p)}_{q,p}(\Omega)=(L_q(\Omega),W^2_q(\Omega))_{1-1/p,p}=(W^2_q(\Omega),L_q(\Omega))_{1/p,p}.    
\end{equation}
We shall not distinguish between notation of spaces for scalar and vector valued functions, i.e. we write $L_q(\Omega)$ instead of $L_q(\Omega)^3$  etc. However, we write vector valued functions in boldface.

We also introduce a brief notation for the regularity class of the solution. Namely, for a triple $(g,h,\vf)$ we define a norm  
\begin{equation} \label{def:X}
\|(g,h,\vf)\|_{{\cal{X}}(T)}=\|\vf\|_{L_p(0,T;W^2_q(\Omega))}+ \| \vf_t\|_{L_p(0,T;L_q(\Omega))}+\|g,h\|_{W^1_p(0,T;W^1_q(\Omega))}
\end{equation}
and a seminorm
\begin{equation} \label{def:dotX}
\begin{aligned}
\|(g,h,\vf)\|_{\dot{\cal{X}}(T)}=&\|\vf\|_{L_p(0,T;W^2_q(\Omega))}+ \| \vf_t\|_{L_p(0,T;L_q(\Omega))}+\|\nabla g,\nabla h\|_{L_p(0,T;L_q(\Omega))}\\
&+\|\d_t g,\d_t h\|_{L_p(0,T;W^1_q(\Omega))}.
\end{aligned}
\end{equation}
Obviously, we denote by ${\cal X}(T)$ a space of functions for which the norm \eqref{def:X} is finite. 

%
Next, for $0<\ep<\pi$ and $\lambda_0>0$ we introduce 
\begin{equation} \label{def:sec}
\begin{aligned}
&\Sigma_\ep = \{ \lambda \in \C \setminus \{0\}: \; |{\rm arg}\lambda|\leq \pi-\ep\}, \quad \Sigma_{\ep,\lambda_0}=\{ \lambda \in \Sigma_\ep: \; |\lambda|\geq\lambda_0 \}, \\ 
&\C_+ := \{ \lambda \in \C: \; {\rm Re}\lambda \geq 0\}.
\end{aligned}
\end{equation}
Finally, by $E(\cdot)$ we shall denote a non-negative non-decreasing continuous function such that $E(0)=0$.  We use the notation $E(T)$ in particular to denote these constants in various inequalities, which can be made arbitrarily small by taking $T$ sufficiently small.

\subsection{Main results}  
The first main result of this paper gives local in time existence and uniqueness for problem \eqref{S}:
\begin{thm} \label{t:local}
Assume that $\Omega$ is a uniform $C^2$ domain and $2<p<\infty, \; 3<q<\infty, \; \frac{2}{p}+\frac{3}{q}<1$.
Assume moreover that $R_0,Q_0$ satisfy \eqref{init:pos}
and $\vu_0$ satisfies the compatibility condition
\begin{equation} \label{comp}
\vu_0|_{\de \Omega}=0.
\end{equation}
Then for any $L>0$ there exists $T>0$ such that if 
\begin{equation} \label{init:loc}
\|\nabla R_0\|_{L_q(\Omega)}+\|\nabla Q_0\|_{L_q(\Omega)}+\|\vu_0\|_{B^{2-2/p}_{q,p}(\Omega)} \leq L
\end{equation}
then \eqref{S} admits a unique solution $(R,Q,\vu)$ on $(0,T)$ with the estimate  
\begin{equation} \label{est:loc}
\|(R,Q,\vu)\|_{\CX(T)} \leq CL, \qquad \int_0^T \|\nabla \vu\|_{L_\infty(\Omega)} \leq \delta,    
\end{equation}
where $\|\cdot\|_{\CX(T)}$ is defined in \eqref{def:X} and $\delta$ is a small positive constant.  
\end{thm}
Next we show global well-posedness for \eqref{S} assuming additionally that $\Omega$ is bounded and initial data is close to some constants: 
\begin{thm} \label{t:global}
Assume that $\Omega$ is a bounded $C^2$ domain, 
$p,q,R_0,Q_0$ satisfy the assumptions of Theorem \ref{t:local}
and let $R_*,Q_*$ be any positive constants. Assume moreover that $\vu_0$ satisfies the compatibility condition \eqref{comp}. Then there exists $\ep>0$ s.t. if 
\begin{equation} \label{small}
\|R_0-R_*\|_{W^1_q(\Omega)}+\|Q_0-Q_*\|_{W^1_q(\Omega)}+\|\vu_0\|_{B^{2-2/p}_{q,p}(\Omega)} \leq \ep,    
\end{equation}
then the solution to \eqref{S} exists globally in time and satisfies the following decay estimate:
\begin{equation} \label{est:glob}
\|e^{\beta t}(R,Q,\vu)\|_{\dot\CX(+\infty)} \leq C\ep, \qquad \int_0^{+\infty} \|\nabla \vu\|_{L_\infty(\Omega)} \leq \delta ,   
\end{equation}
where $\beta$ is a positive constant, $\delta$ is a small positive constant, and $\|\cdot\|_{\dot\CX}$ is defined in \eqref{def:dotX}. 
\end{thm}

\subsection{Discussion}
The existence of global weak solutions to \eqref{spec} in a space-periodic setting was established by Bresch et al. in \cite{BDGG}, however with additional capillarity effects.
In a one dimensional case the existence of global weak solutions without capillarity was shown in \cite{BHL}, where, moreover, it is shown that after a finite time a vacuum is excluded and at least component corresponding to one phase becomes regular. More recently, the global finite energy weak solutions to a two-fluid Stokes system has been shown by Bresch, Mucha and Zatorska in \cite{BMZ}. Soon after, similar result for full two-fluid compressible Navier-Stokes equations with single velocity was proven by Novotn\'y and  Pokorn\'y in \cite{NoPo}, and by Jin and Novotn\'y \cite{JiNo} for the so-called Baer-Nunziato type of system.

The existence of regular solutions for a simplified version of \eqref{spec_1} where influence of one of the phases in the convective term is neglected has been shown in \cite{Guo}. Existence of global strong solutions in $L_2$ framework in the whole space in $\R^3$ under certain smallness assumptions has been shown by Evje, Wang and Wen in \cite{Evje1}. In the one dimensional case, more analytical and numerical results can be found in \cite{Evje2}, \cite{Evje3}.  Existence, uniqueness and stability of global weak solutions, still in one dimension, to system  \eqref{spec}, has been shown by Li, Sun, Zatorska in \cite{LiSuZa} and some non-uniqeness results for the inviscid version of this system in multiple dimensions was shown by Li and Zatorska in \cite{LiZa}. Finally, let us also mention that very recently Wang, Wen and Yao proved global existence and uniqueness of solutions to the non-conservative two-fluid model with unequal velocities using the energy methods \cite{WaWeYa}.

In order to prove Theorem \ref{t:local} we show a maximal regularity estimate for associated linear problem. For this purpose we apply the theory based on the concept of $\CR$-boundedness. The underlying result for this approach is the famous Weis' vector-valued Fourier  multiplier theorem proved in \cite{Weis}, which we recall in the Appendix. In brief, this result allows to deduce maximal regularity for a time-dependent problem from $\CR$-boundedness of a family of solution operators to associated resolvent problem. In \cite{DHP} this approach was applied to prove maximal regularity for the heat equation using explicit solution formula. Further development of the theory in the context of equations of fluid mechanics was mostly due to Y. Shibata and his collaborators. In \cite{SS2} the authors show a maximal regularity for the Stokes problem with Neumann boundary condition. For this purpose they prove several technical results enabling to show $\CR$ boundedness for the resolvent problem. These results have been applied and extended in \cite{ES} to treat the compressible Navier-Stokes system with Dirichlet boundary condition.
The result has been extended to slip boundary conditions in \cite{Murata}, \cite{MS16}.
For further developments we can mention \cite{S1},\cite{SS1} for free boundary problems, \cite{KSK} for a 2-phase compressible-incompressible flow or \cite{HM} for a compressible fluid-rigid body interaction problem.
In a series of papers of Piasecki, Shibata and Zatorska \cite{PSZ}-\cite{PSZ3} the above described approach was applied to treat a system describing multi-component compressible mixtures. Here we rely on the ideas we developed in these papers to treat quite complex nonlinearities resulting from Lagrangian transformation and linearization.

\section{Lagrangian coordinates}
In order to define the Lagrangian transformation we start with a following simple observation:
\begin{lem} \label{l:lag1}
Let $p$ and $q$ satisfy the assumptions of Theorem \ref{t:local}. Then\\
(i) if  
$\|f\|_{L_p(0,T;W^2_q(\Omega))} \leq M$
for some $M>0$, then
\begin{equation} \label{2:1}
\int_0^T \|\nabla f(t,\cdot)\|_{L_\infty(\Omega)} \dt \leq M E(T), \end{equation}
\noindent
(ii) if 
$\|e^{\gamma t}f\|_{L_p(0,\infty;W^2_q(\Omega))}\leq M$ 
for some $M,\gamma>0$ then
\begin{equation} \label{2:1a}
\int_0^\infty \|\nabla f(t,\cdot)\|_{L_\infty(\Omega)} \dt \leq CM. \end{equation}
\end{lem}
{\em Proof.}
By the imbedding theorem and H\"older inequality we have 
\begin{align*}
\int_0^T\|\nabla f(t,\cdot)\|_{L_\infty(\Omega)} \dt \leq C \int_0^T \|f(t,\cdot)\|_{W^2_q(\Omega)}\dt \leq T^{1/p'}\int_0^T \left(\|f(t,\cdot)\|_{W^2_q(\Omega)}^p \right)^{1/p} \dt \leq ME(T),   
\end{align*}
which proves the first assertion, and for the second we have 
\begin{align*}
\int_0^\infty\|\nabla f(t,\cdot)\|_{L_\infty(\Omega)} \dt\leq 
C \int_0^\infty e^{-\gamma t}e^{\gamma t}\|f(t,\cdot)\|_{W^2_q(\Omega)} \dt\\ \nonumber
\leq \left( \int_0^\infty e^{-\gamma tp'} \dt \right)^{1/p'}
\|e^{\gamma t}f\|_{L_p(0,\infty;W^2_q(\Omega))}\leq CM.
\end{align*}
\qed

Now we can proceed with the transformation. Let $\vv(y, t)$ be the 
velocity field in the Lagrange coordinates, so that the following change of variables is true:
\begin{equation}\label{lag:1}
x = y + \int^t_0\vv(y, s)\,ds.
\end{equation}
Then for any differentiable function $f$ we have 
\begin{equation} \label{dt_lag}
\pt f(t,\phi(t,y))=\pt f+\vv \cdot \nabla_x f.
\end{equation}
Moreover, we have
\begin{equation}\label{lag:2}
\frac{\pd x_i}{\pd y_j} = \delta_{ij} 
+ \int^t_0\frac{\pd v_i}{\pd y_j}(y, s)\,ds,
\end{equation}
where $\delta_{ij}$ are Kronecker's delta symbols. 
If $\vv$ is the solution satifying the regularity from Theorem \ref{t:local}, then by Lemma \ref{l:lag1} we can assume
\begin{equation}\label{assump:1}
\sup_{t \in (0,T)}\int^t_0\|\nabla\vv(\cdot, s)\|_{L_\infty(\Omega)}\,ds
\leq \delta
\end{equation}
with some small positive constant $\delta$. Therefore by \eqref{lag:2} the $N\times N$ matrix
$\pd x/\pd y = (\pd x_i/\pd y_j)$ has the inverse 
\begin{equation}\label{lag:3}
\lr{\frac{\pd x_i}{\pd y_j}}^{-1} = \bI + \bV^0(\vk_\vv)
\end{equation}
where 
\begin{equation} \label{kv}
\vk_\vv = \int^t_0\nabla\vv(y, s)\,ds,    
\end{equation} 
$\bI$ is the $3\times 3$ identity matrix, and $\bV^0(\vk)$ is the $3\times 3$ matrix of 
smooth functions with respect to $\vk \in \R^{3\times3}$ defined on $|\vk| <\delta$ with $\bV^0(0) = 0$,
where $\vk$ are independent variables corresponding to 
$\vk_\vv$. 

We have
\begin{equation}\label{lag:4}
\nabla_x = (\bI + \bV^0(\vk_\vv))\nabla_y, 
\quad \frac{\pd}{\pd x_i} = \sum_{j=1}^N (\delta_{ij} + V^0_{ij}(\vk_\vv))
\frac{\pd}{\pd y_j}.
\end{equation}
Notice that adding \eqref{SR} to \eqref{ST} and using this equation to compute $\de_t(R+Q)$ we can rewrite \eqref{SM} as 
\begin{equation}
(R+Q)\d_t\vu+(R+Q)\vu \cdot \nabla \vu-\Div \vS(\vu) +\Grad Z^{\gamma^+}=\vc{0}.    
\end{equation}
Moreover, the map:
$x = \Phi(y, t)$ is bijection from $\Omega$ onto
$\Omega$, and so setting
\begin{equation}\label{lag:5}
\vv(y, t) = \vu(x, t),
\quad \fr(y, t) = R(x, t), \quad
\fq(y, t) = Q(x, t), \quad \fz(y,t)=Z(x,t),
\end{equation}
we can transform \eqref{S} to the following form
\begin{subequations} \label{sysL1}
\begin{align}
&\pt \fr+\fr\Div\vv=O_1(U) \label{sysL1:1}\\[3pt]
&\pt \fq+\fq\Div\vv=O_2(U) \label{sysL1:2}\\[3pt]
&(\fr+\fq)\pt\vv-\mu\lap\vv-\nu\Grad\Div\vv+ \gamma^+\Big(\fz^{\gamma^+-1}\partial_\fr \fz\Grad \fr+ \fz^{\gamma^+-1}\partial_ \fq \fz\Grad\fq\Big)=\vc{O}_3(U)\label{sysL1:3},\\[3pt]
& \fq=\lr{1-\frac{ \fr}{ \fz}} \fz^\gamma,\label{sysL1:4}\\[3pt]
& \fr|_{t=0}=\fr_0, \quad  \fq|_{t=0}=\fq_0, \quad \vv|_{t=0}=\vv_0,\\[3pt]
&\vv|_{\d\Omega}=0
\end{align}
\end{subequations}
where 
\begin{equation} \label{def:U}
U=(\vv,\fr, \fq).    
\end{equation}
For consistency with notation \eqref{lag:5} let us denote the initial data for the problem in Lagrangian coordinates as  
$$
\fr_0(y)=R_0(y), \quad \fq_0(y)=Q_0(y),\quad \vv_0(y)=\vu_0(0). 
$$
We derive explicit form of the terms $\de_ \fr \fz, \, \de_ \fq \fz$.
Recall that $ \fr \leq  \fz$. Let us now apply $\partial_{\fq}$, $\partial_{\fr}$ to both sides of \eqref{sysL1:4}, we obtain respectively:
\[1=\gamma  \fz^{\gamma-1}\partial_{\fq}  \fz- \fr(\gamma-1) \fz^{\gamma-2}\partial_{\fq}  \fz,\]
\[0=\gamma  \fz^{\gamma-1}\partial_{\fr}  \fz- \fz^{\gamma-1}- \fr(\gamma-1) \fz^{\gamma-2}\partial_{\fr}  \fz,
\]
therefore,
\eq{\label{pTR}
\partial_{\fq}  \fz=\frac{1}{\gamma  \fz^{\gamma-1}- \fr(\gamma-1) \fz^{\gamma-2}},\quad
\partial_{\fr}  \fz= \frac{ \fz^{\gamma-1}}{\gamma  \fz^{\gamma-1}- \fr(\gamma-1) \fz^{\gamma-2}},
}
and so 
\begin{equation} \label{pTR:1}
 \fz^{\gamma^+-1}\partial_{\fq} \fz=\frac{\fz^{\gamma^+}}{\gamma \fz^\gamma-(\gamma-1) \fr\fz^{\gamma-1}}, \qquad
 \fz^{\gamma^+-1}\partial_{\fr} \fz=\frac{\fz^{\gamma^+}}{\gamma \fz-(\gamma-1) \fr}.    
\end{equation}
Due to \eqref{a:gamma} and \eqref{RleqZ} we have
\begin{equation} \label{zeta:1}
\frac{1}{\gamma  \fz^{\gamma-1}} \; \leq \; \d_ \fq  \fz \; \leq \; \frac{1}{\fz^{\gamma-1}}, \qquad \gamma^{-1}\;\leq \; \d_\fr  \fz \; \leq \; 1    
\end{equation}
and 
\begin{equation}
\gamma^{-1} \fz^{\gamma^+ -\gamma} \; \leq \;  \fz^{\gamma^+-1}\d_ \fq  \fz \; \leq \;  \fz^{\gamma^+ -\gamma}, \qquad 
\gamma^{-1}  \fz^{\gamma^+ -1} \; \leq \;  \fz^{\gamma^+-1}\d_ \fr  \fz \; \leq \;  \fz^{\gamma^+ -1}.
\end{equation}
We have  
\begin{equation}\label{lag:div}
\Div_x = \Div_y + \sum_{i,j=1}^3 V^0_{ij}(\vk_\vv)\frac{\pd v_i}{\pd y_j},
\end{equation} 
therefore by \eqref{dt_lag},\eqref{lag:4} and \eqref{lag:5}, we obtain \eqref{sysL1:1} with
\begin{equation}\label{lag:6}
O_1(U) = - \fr\sum_{i,j=1}^3 V^0_{ij}(\vk_\vv)\frac{\pd v_i}{\pd y_j}.
\end{equation}
Similarly,
\begin{equation}\label{lag:6b}
O_2(U) = - \fq\sum_{i,j=1}^3 V^0_{ij}(\vk_\vv)\frac{\pd v_i}{\pd y_j}.
\end{equation}
To find $\vc{O}_3(U)$ we need to transform the second order operators.
By \eqref{lag:4}, we have
$$
\Delta \vu = \sum_{k=1}^3\frac{\pd}{\pd x_k}\lr{\frac{\pd \vu}{\pd x_k}}
= \sum_{k,\ell,m=1}^3\lr{\delta_{k\ell} + V^0_{kl}(\vk_\vv)}
\frac{\pd}{\pd y_\ell}
\lr{(\delta_{km} + V^0_{km}(\vk_\vv))\frac{\pd \vv}{\pd y_m}},
$$
and so 
setting 
\eqh{
A_{2\Delta}(\vk)\nabla^2\vv &= 2\sum_{\ell, m=1}^3 V^0_{k\ell}(\vk)
\frac{\pd^2\bv}{\pd y_\ell\pd y_m}
+ \sum_{k,\ell, m=1}^3
V^0_{k\ell}(\vk)V^0_{km}(\vk)
\frac{\pd^2\vv}{\pd y_\ell \pd y_m}, \\
A_{1\Delta}(\vk)\nabla\vv & = \sum_{\ell, m=1}^3(\nabla_\vk V^0_{\ell m})(\vk)
\int^t_0(\pd_l\nabla\vv)\,ds \frac{\pd \vv}{\pd y_m}
+ \sum_{k, \ell, m=1}^3
V^0_{k\ell}(\vk) (\nabla_\vk V^0_{km})(\vk)
\int^t_0\pd_\ell\nabla\vv\,ds\frac{\pd \vv}{\pd y_m}
}
we have
$$\Delta \vu = \Delta \vv + A_{2\Delta}(\vk_\vv)\nabla^2\vv
+ A_{1\Delta}(\vk_\vv)\nabla\vv.
$$
Moreover, by \eqref{lag:4}, we have
$$\frac{\pd}{\pd x_j}\Div\vu 
= \sum_{k=1}^3(\delta_{jk} + V^0_{jk}(\vk_\vv))\frac{\pd}{\pd y_k}
\lr{\Div\vv + \sum_{\ell, m=1}^3 V^0_{\ell m}(\vk_\vv)\frac{\pd v_\ell}{\pd y_m}},
$$
therefore setting
\eqh{
A_{2\Div,j}(\vk)\nabla^2\vv
& = \sum_{\ell, m=1}^3V^0_{\ell m}(\vk)\frac{\pd^2 v_\ell}{\pd y_m\pd y_j}
+ \sum_{k=1}^3 V^0_{jk}(\vk)\frac{\pd}{\pd y_k}\Div\vv
+ \sum_{k, \ell=1}^3V^0_{jk}(\vk)V^0_{\ell m}(\vk)
\frac{\pd^2v_\ell}{\pd y_k\pd y_m}, \\
A_{1\Div, j}(\vk)\nabla\bv
& = \sum_{\ell, m=1}^3(\nabla_\vk V^0_{\ell m})(\vk)
\int^t_0\pd_j\nabla\vv\,ds\frac{\pd v_\ell}{\pd y_m} 
+ \sum_{k,\ell, m=1}^3V^0_{jk}(\vk)(\nabla_\vk V^0_{\ell m})(\vk)
\int^t_0\pd_k\nabla\vv\,ds\frac{\pd v_\ell}{\pd y_m},
}
we obtain
$$\frac{\pd}{\pd x_j}\Div\vu
= \frac{\pd}{\pd y_j}\Div\vv + A_{2\Div, j}(\vk_\vv)\nabla^2\vv
+ A_{1\Div, j}(\vk_\vv)\nabla\vv.
$$
Since
\eqh{
&Z^{\gamma^+-1}\partial_RZ\Grad R+Z^{\gamma^+-1}\partial_QZ\Grad Q\\
&= \fz^{\gamma^+-1}\partial_ \fr \fz(\Grad \fr+ \bV^0(\vk_\vv)\nabla \fr)+ \fz^{\gamma^+-1}\partial_\fq \fz(\Grad \fq+ \bV^0(\vk_\vv)\nabla \fq),
}
we have
\begin{equation}\label{lag:7}\begin{split}
 \vc{O}_3(U) &= \mu A_{2\Delta}(\vk_\vv)\nabla^2\vv 
+ \mu A_{1\Delta}(\vk_\vv)\nabla\vv
+ \nu A_{2\Div}(\vk_\vv)\nabla^2\vv + \nu A_{1\Div}(\vk_\vv)\nabla\vv \\
&-  \fz^{\gamma^+-1}\partial_\fr \fz \bV^0(\vk_\vv)\nabla \fr
- \fz^{\gamma^+-1}\partial_\fq \fz \bV^0(\vk_\vv)\nabla \fq.
\end{split}\end{equation}

\section{Local well-posedness}
\subsection{Linearization around the initial condition}

We introduce new unknowns 
\begin{equation} \label{pert:1}
\sigma = \fr-\fr_0, \quad \eta = \fq-\fq_0.
\end{equation}
Denoting 
\begin{equation}
V=(\vv,\sigma,\eta)=U-(0,\fr_0,\fq_0)    
\end{equation}
we obtain the following system with time-independent coefficients
\eq{\label{sysL2}
&\pt\sigma+\fr_0\Div\vv=f_1(V),\\
&\pt\eta+\fq_0\Div\vv=f_2(V),\\
&(\fr_0+\fq_0)\pt\vv-\mu\lap\vv-\nu\Grad\Div\vv+\omega_1^0\Grad\sigma+\omega_2^0\Grad\eta=\vf_3(V),\\
&(\sigma,\eta,\vv)|_{t=0}=(0,0,\vv_0), \quad \vv|_{\d\Omega}=\vc{0},
}
where $ \fz_0= \fz_0(\fr_0,\fq_0)$ is defined by 
\begin{equation*} 
\fq_0=\left(1-\frac{\fr_0}{ \fz_0} \right) \fz_0^\gamma    
\end{equation*}
and
\begin{equation} \label{gamma12:1} 
\omega_1^0 = \frac{ \gamma^+\fz_0^{\gamma^+}}{\gamma \fz_0-(\gamma-1)\fr_0}, \qquad
\omega_2^0 = \frac{ \gamma^+\fz_0^{\gamma^+}}{\gamma \fz_0^{\gamma}-(\gamma-1)\fr_0 \fz_0^{\gamma-1}}.
\end{equation}
The right hand side of \eqref{sysL2} is given by 
\begin{align}
f_1(V)&=O_1(V+(0,\fr_0,\fq_0))-\sigma\Div\vv, \label{f1} \\[5pt]
f_2(V)&=O_2(V+(0,\fr_0,\fq_0))-\eta\Div\vv, \label{f2} \\[8pt]
\vf_3(V)&=\vc{O}_3(V+(0,\fr_0,\fq_0))-(\sigma+\eta)\d_t\vv \label{f3}\\[5pt]\nonumber
&-\underbrace{ \frac{ \gamma^+\gamma\big[ \fz_0(  \fz^{\gamma^+}- \fz_0^{\gamma^+} )+ \fz_0^{\gamma^+}( \fz_0- \fz)\big]+(\gamma-1)\big[  \fz_0^{\gamma^+}( \fr- \fr_0)+ \fr_0( \fz_0^{\gamma^+}- \fz^{\gamma^+}) \big] }{\big(\gamma \fz_0-(\gamma-1) \fr_0\big)\big(\gamma \fz-(\gamma-1) \fr\big) }\nabla \sigma }_{I_1( \fr, \fq,\fz)} \\[5pt]
\nonumber& - \underbrace{ \frac{ \gamma^+\gamma\big[  \fz_0^{\gamma^+}( \fz_0^{\gamma}- \fz^{\gamma})+ \fz_0^{\gamma}( \fz^{\gamma^+}- \fz_0^{\gamma^+}) \big] }{ \big(\gamma \fz^{\gamma}-(\gamma-1) \fr \fz^{\gamma-1}\big)\big(\gamma \fz_0^{\gamma}-(\gamma-1) \fr_0 \fz_0^{\gamma-1}\big)}\nabla\eta}_{I_2(\fr, \fq,\fz)}\\[5pt]
\nonumber& -\underbrace{\frac{\gamma^+(\gamma-1)\big[ \fr \fz_0^{\gamma^+}(\fz^{\gamma-1}- \fz_0^{\gamma-1})+ \fz_0^{\gamma^+} \fz_0^{\gamma-1}(\fr- \fr_0)+ \fr_0 \fz_0^{\gamma-1}( \fz_0^{\gamma^+}- \fz^{\gamma^+}) \big]  }{ \big(\gamma \fz^{\gamma}-(\gamma-1) \fr \fz^{\gamma-1}\big)\big(\gamma \fz_0^{\gamma}-(\gamma-1) \fr_0 \fz_0^{\gamma-1}\big)}\nabla\eta}_{I_3( \fr, \fq, \fz)}  \\[5pt]
\nonumber& -\underbrace{\left(\frac{ \gamma^+\fz^{\gamma^+}}{\gamma \fz-(\gamma-1) \fr}\nabla \fr_0
+\frac{ \gamma^+\fz^{\gamma^+}}{\gamma \fz^{\gamma}-(\gamma-1) \fr \fz^{\gamma-1}}\nabla \fq_0\right)}_{I_4(\fr, \fq, \fz)}  .
\end{align}

\subsection{Maximal regularity}
In this section we show the maximal regularity for the linear problem corresponding to \eqref{sysL2}, i.e.:
\eq{\label{sys:lin1}
&\pt\sigma+\fr_0\Div\vv=f_1,\\
&\pt\eta+\fq_0\Div\vv=f_2,\\
&(\fr_0+\fq_0)\pt\vv-\mu\lap\vv-\nu\Grad\Div\vv+\omega_1^0\Grad\sigma+\omega_2^0\Grad\eta=\vf_3,\\
&(\sigma,\eta,\vv)|_{t=0}=(\sigma_0,\eta_0,\vv_0), \quad \vv|_{\d\Omega}=0.
}
The main result of this section reads\\
\begin{thm} \label{t:maxreg}
Let $\Omega\subset \R^3$ be a uniform $C^2$ domain.
Let $1<p,q<\infty$, $\frac{2}{p}+\frac{1}{q} \neq 1$.
Assume that there exists constants 
$a_1,a_2,\kappa>0$ and $r>3$ such that
\begin{equation} \label{reg:init}
\begin{aligned}
0 \leq \fr_0, \fq_0, \omega_1^0, \omega_2^0 \leq a_1, \quad \fr_0+\fq_0\geq \kappa,\\[3pt]
\|\nabla \fr_0, \nabla \fq_0, \nabla \omega_1^0, \nabla \omega_2^0\|_{L_r(\Omega)} \leq a_2.
\end{aligned}
\end{equation}

Assume moreover that
\begin{align*}
f_1,f_2 \in L_p(0,T,W^1_q(\Omega)),\; \vf_3 \in L_p(0,T,L_q(\Omega)),\;
\sigma_0,\eta_0 \in W^1_q(\Omega), \; \vv_0 \in B^{2-1/p}_{q,p}(\Omega). 
\end{align*}
Finally, for $\frac{2}{p}+\frac{1}{q}<2$ assume that $\vv_0$ satisfy the compatibility condition \eqref{comp}.
Then there exists positive constants $b,C$ such that for any $T>0$ problem \eqref{sys:lin1} admits a unique solution $(\sigma,\eta,\vv) \in {\cal X}(T)$ with the estimate 
\eq{ \label{est:linME1_A}
\|(\sigma,\eta,\vv)\|_{{\cal X}(T)}\leq C e^{b T}&[\|\vv_0\|_{B^{2-2/p}_{q,p}(\Omega)}+\|\sigma_0,\eta_0\|_{W^1_q(\Omega)}\\[3pt]
&+\|\vf_3\|_{L_p(0,T;L_q(\Omega))}+\|(f_1,f_2)\|_{L_p(0,T;W^1_q(\Omega))}].
}
\end{thm}

In order to prove Theorem \ref{t:maxreg} we prove $\CR -boundedness$ for associated resolvent problem which is obtained applying Laplace transform to \eqref{sys:lin1}: 
\begin{subequations} \label{sys:res}
\begin{align}
&\lambda \sigma_\lambda+ \fr_0\Div\vv_\lambda=f_1,\label{res1.1}\\
&\lambda \eta_\lambda+ \fq_0\Div\vv_\lambda=f_2,\label{res1.2}\\
&( \fr_0+ \fq_0)\lambda\vv_\lambda-\mu\lap\vv_\lambda-\nu\Grad\Div\vv_\lambda+\omega_1^0\Grad\sigma_\lambda+\omega_2^0\Grad\eta_\lambda=\vf_3,\label{res1.3}\\
&\vv_\lambda|_{\d\Omega}=0.
\end{align}
\end{subequations}
In the Appendix we recall how $\CR -boundedness$ for the resolvent problem implies maximal regularity due to Weis vector-valued Fourier Multiplier Theorem. Therefore in order to prove Theorem \ref{t:maxreg} it is sufficient to show the following result: 
\begin{thm} \label{t:rbound}
Let $1 < q < \infty$ and $0 < \epsilon < \pi/2$.
Assume that $\Omega,\fr_0,\fq_0,\omega_1^0,\omega_2^0$ satisfy the assumptions of Theorem \ref{t:maxreg}. 
Then, there exist a positive constant $\lambda_0$ and 
operator families $\CA_1(\lambda),\CA_2(\lambda) \in {\rm Hol}\,(\Lambda_{\epsilon, \lambda_0},
\CL(W^1_q(\Omega)^2 \times L_q(\Omega), W^1_q(\Omega)))$ and  
$\CB(\lambda) \in {\rm Hol}\,(\Lambda_{\epsilon, \lambda_0},
\CL(W^1_q(\Omega)^2 \times L_q(\Omega), W^2_q(\Omega)^N))$, 
such that for any $(f_1,f_2,\vf_3) \in  (W^1_q(\Omega))^2\times L_q(\Omega)$
and $\lambda \in \Sigma_{\epsilon, \lambda_0}$, 
$(\sigma_\lambda = \CA_1(\lambda)(f_1,f_2,\vf_3)$, $\eta_\lambda = \CA_2(\lambda)(f_1,f_2,\vf_3)$, $ 
\vu_\lambda= \CB(\lambda)(f_1,f_2,\vf_3))$
is a unique solution of \eqref{sys:res} and 
\begin{align*}
\CR_{\CL( (W^1_q(\Omega))^2 \times L_q(\Omega), W^1_q(\Omega))}(\{(\tau\pd_\tau)^\ell
\CA_1(\lambda):\; \lambda \in \Sigma_{\epsilon, \lambda_0}\}) &\leq M, \\
\CR_{\CL( (W^1_q(\Omega))^2 \times L_q(\Omega), W^1_q(\Omega))}(\{(\tau\pd_\tau)^\ell
\CA_2(\lambda):\; \lambda \in \Sigma_{\epsilon, \lambda_0}\}) &\leq M, \\
\CR_{\CL(L_q(\Omega) \times W^1_q(\Omega), W^{2-j}_q(\Omega)^N)}(\{(\tau\pd_\tau)^\ell
(\lambda^{j/2}\CB(\lambda)):\; 
\lambda \in \Sigma_{\epsilon, \lambda_0}\}) &
\leq M,
\end{align*}
for $\ell=0,1$, $j=0,1,2$ and some constant $M>0$, where { by ${\rm Hol}$ we denote the space of holomorphic operators.}
\end{thm}

\begin{rmk}
We say that the operator valued function $T(\lambda)$ is holomorphic if it is differentiable in norm for all $\lambda$ in a complex domain. For more details, see \cite[Chaper VII,\$1.1]{Kato}.
\end{rmk}
Following the idea introduced already in \cite{ES} in context of the compressible Navier-Stokes equations, in order to prove Theorem \ref{t:rbound} we transform \eqref{sys:res} using the fact that continuity equation becomes an algebraic equation in the resolvent problem. Computing $\eta_\lambda$ and $\sigma_\lambda$ from \eqref{res1.1}-\eqref{res1.2} we obtain 
\begin{equation*}
\sigma_\lambda = \lambda^{-1}\left(f_1- \fr_0\div\vv\right), \quad 
\eta_\lambda = \lambda^{-1}\left(f_2- \fq_0\div\vv\right).
\end{equation*}
Plugging these indentities into \eqref{res1.3} we obtain 
\begin{equation} \label{res:ME}
\begin{aligned}
&( \fr_0+ \fq_0)\lambda\vv_\lambda-\mu\Delta\vv_\lambda
-[ \nu+\lambda^{-1}( \omega_1^0 \fr_0+\omega_2^0 \fq_0 ) ]\nabla \div\vv\\
&=\vf_3+\lambda^{-1}\omega_1^0(\div\vv\nabla \fr_0-\nabla f_1)
+\lambda^{-1}\omega_2^0(\div\vv\nabla \fq_0-\nabla f_2).
\end{aligned}
\end{equation}
Since the RHS of \eqref{res:ME} contains lower order terms, it is enough to 
show $\CR$-boundedness for a problem 
\begin{equation} \label{res:ME2}
( \fr_0+ \fq_0)\lambda\vv_\lambda-\mu\Delta\vv_\lambda
-[ \nu+\lambda^{-1}( \omega_1^0 \fr_0+\omega_2^0 \fq_0 ) ]\nabla \div\vv =\vf.  \end{equation}
The result is 
\begin{prop} \label{prop:Rbound} Let $1 < q < \infty$ and 
$0 < \epsilon < \pi/2$.  Assume that $\Omega$ is a uniform 
$C^2$ domain in $\R^N$.  Then, there exists a positive constant
$\lambda_0$ such that there exists   
an operator family $\CC(\lambda) \in {\rm Hol}\,
(\Sigma_{\epsilon, \lambda_0}, \CL(L_q(\Omega), W^2_q(\Omega)))
$
such that for any 
$\lambda \in \Sigma_{\epsilon, \lambda_0}$
and $\vf \in L_q(\Omega)$, 
$\bv = \CC(\lambda)\vf$ is a unique 
solution of \eqref{res:ME2}, and
$$\CR_{\CL(L_q(\Omega), W^{2-j}_q(\Omega))}
(\{(\tau\pd_\tau)^\ell\CC(\lambda) \mid \lambda \in 
\Sigma_{\epsilon, \lambda_0}\}) \leq M$$
for $\ell=0,1$, $j=0,1,2$ and some constant $M>0$.
\end{prop}
{\bf Proof.} An analog of Proposition \ref{prop:Rbound} has been shown in (\cite{ES}, Theorem 2.10)  for a problem 
\begin{equation} \label{ME:resES}
\lambda \vu - \mu \Delta_y \vu - (\nu+\gamma^2\lambda^{-1}) \nabla_y \div_y \vu  = \vf , \quad
\vu|_{\de \Omega}=0    
\end{equation}
where $\mu,\nu$ and $\gamma$ are constants satisfying $\mu+\nu>0$ and 
$\gamma>0$. The proof requires only minor modifications in order to prove Proposition \ref{prop:Rbound}, therefore we present only a sketch. First we solve a problem with constant coefficients in the whole space
\begin{equation} \label{ME:res:const}
( \fr_0^*+ \fq_0^*)\lambda \vv_\lambda-\mu\Delta\vv_\lambda
-[ \nu+\lambda^{-1}( \omega_1^{0*} \fr_0^*+\omega_2^{0*} \fq_0^* ) ]\nabla \div\vv =\vf \quad \in \R^n,
\end{equation}
where $ \fr_0^*, \fq_0^*,\omega_1^{0*},\omega_2^{0*}$ are constants satisfying 
\begin{equation*}
\fr_0^*, \,  \fq_0^*, \omega_1^{0*},\omega_2^{0*} \geq 0, \quad  \fr_0^* +  \fq_0^* > 0.
\end{equation*}
$\CR$-boundedness for \eqref{ME:res:const} can be shown following the proof Theorem 3.1 in \cite{ES}. The latter is shown for problem \eqref{ME:resES} in the whole space. In order to adapt it to \eqref{ME:res:const} we can divide \eqref{ME:res:const} by $ \fr_0^*+ \fq_0^*$ using the fact that this constant is strictly positive. We obtain  
\begin{equation*}
\lambda\vv_\lambda-\mu^*\Delta\vv_\lambda
-[ \nu^*+\lambda^{-1}\omega^* ]\nabla \div\vv =\vf \quad \in \R^n,    
\end{equation*}
where
$$
\mu^*= \frac{\mu}{ \fr_0^*+ \fq_0^*}, \quad \nu^*= \frac{\nu}{ \fr_0^*+ \fq_0^*}, \quad  
\omega^*=\frac{ \omega_1^{0*} \fr_0^*+\omega_2^{0*} \fq_0^*}{ \fr_0^*+ \fq_0^*}. 
$$
Therefore the only difference is that now we have $\omega^* \geq 0$ instead of $\gamma^2$. However, strict positivity of this constant neither its square structure is not necessary, nonnegativity is sufficient. 

Next we consider \eqref{ME:res:const} in a half-space supplied with the boundary condition $\vv|_{\de \R^n_+}=0$. Here we can follow the proof of Theorem 4.1 in \cite{ES} which works without modifications for $\omega^* \geq 0$. 

The third step consists in showing $\CR$-boundedness in a bent half-space, the necessary result is Theorem 5.1 in \cite{ES}, where replacing $\lambda^{-1}\gamma^2>0$ with  
$\omega^* \geq 0$ is again harmless. 

The final step is application of partition of unity and properties of a uniform $C^2$ domain. Here we have to deal with variable coefficients which is not in the scope of Theorem 2.10 in \cite{ES}. 
However, we can refer to a more recent result, Theorem 4.1 in \cite{PSZ} which gives $\CR$-boundedness for a resolvent problem corresponding to more complicated system describing flow of a two-component mixture. It is sufficient to follow  Section 6.3 of \cite{PSZ} with some obvious modifications.
At this stage we need continuity of the coefficients $\fr_0,\fq_0,\omega_1^0,\omega_2^0$ which is assured by \eqref{reg:init}. 
\qed

\subsection{Preliminary estimates}
We start with recalling two embedding results for Besov spaces.  
The first one is \cite[Theorem 7.34 (c)]{Ad}:
\begin{lem} \label{l:imbed1}
Assume $\Omega \subset \R^n$ satisfies the cone condition and let $1\leq p,q\leq\infty$ and $sq>n$. Then 
$$
B^s_{q,p}(\Omega) \subset C_B(\Omega),
$$
where $C_B$ we denote the space of continuous bounded functions.
\end{lem}
In particular $u \in B^{2-2/p}_{q,p}(\Omega)$ implies $\nabla u \in B^{1-2/p}_{q,p}(\Omega)$. Therefore Lemma \ref{l:imbed1} with $s=1-2/p$ and $n=3$ yields 
\begin{cor} \label{c:imbed1}
Assume $\frac{2}{p}+\frac{3}{q}<1$ and let $\Omega$ satisfy the assumptions of Theorem \ref{t:local}. Then $B^{2-2/p}_{q,p}(\Omega) \subset W^1_{\infty}(\Omega)$ and 
\begin{equation}
\|f\|_{W^1_\infty(\Omega)} \leq C \|f\|_{B^{2-2/p}_{q,p}(\Omega)}.    
\end{equation}
\end{cor}
The next result is due to Tanabe  
(cf. \cite[p.10]{Tanabe}):
\begin{lem} \label{L:int}
Let $X$ and $Y$ be two Banach spaces such that
$X$ is a dense subset of $Y$ and $X\subset Y$ is continuous.
Then for each $p \in (1, \infty)$  
$$W^1_p((0, \infty), Y) \cap L_p((0, \infty), X) 
\subset C([0, \infty), (X, Y)_{1/p,p})$$
and for every $u\in W^1_p((0, \infty), Y) \cap L_p((0, \infty), X)$ we have
$$\sup_{t \in (0, \infty)}\|u(t)\|_{(X, Y)_{1/p,p}}
\leq (\|u\|_{L_p((0, \infty),X)}^p
+ \|u\|_{W^1_p((0, \infty), Y)}^p)^{1/p}.
$$
\end{lem}
This result allows to show the following imbedding.
\begin{lem} \label{l:nonlin1}
Assume $p,q$ satisfy the assumptions of Theorem \ref{t:local}.  
Let $f_t \in L_p(0,T;L_q(\Omega))$,
$f \in L_p(0,T;W^2_q(\Omega))$,  $f(0,\cdot)\in B^{2-2/p}_{q,p}(\Omega)$. Then 
\begin{align} \label{4.1}
\sup_{t \in (0,T)} \|f\|_{B^{2(1-1/p)}_{q,p}(\Omega)} \leq C[ \|f_t\|_{L_p(0,T;L_q(\Omega))}+\|f\|_{L_p(0,T;W^2_q(\Omega))} +\|f(0)\|_{B^{2-2/p}_{q,p}(\Omega)}], \\[5pt] \label{4.2}
\sup_{t \in (0,T)} \|f\|_{W^1_\infty(\Omega))} \leq C[ \|f_t\|_{L_p(0,T;L_q(\Omega))}+\|f\|_{L_p(0,T;W^2_q(\Omega))} +\|f(0)\|_{B^{2-2/p}_{q,p}(\Omega)}],
\end{align}
where $C$ does not depend ot $T$.
\end{lem}
{\bf Proof.} In order to prove \eqref{4.1} we introduce
an extension operator
\begin{equation} \label{def:ext} e_T[f](\cdot, t)
= \begin{cases}
f(\cdot, t) \quad &t\in(0,T), \\
f(\cdot, 2T-t)\quad & t\in(T,2T),\\
0 \quad & t\in(2T,+\infty),
\end{cases}
\end{equation}
If $f|_{t=0}=0$, then we have
\begin{equation} \label{ext:2} \pd_te_T[f](\cdot, t)
= \begin{cases}
0 \quad &t\in(2T,+\infty),\\
(\pd_tf)(\cdot, t) \quad &t\in(0,T), \\
-(\pd_tf)(\cdot, 2T-t)\quad & t\in(T,2T) \\
\end{cases}
\end{equation}
in a weak sense. 
Therefore obviously 
\begin{equation} \label{ef:norm}
\|e_T f\|_{L_p(\R_+;W^2_q(\Omega))}+
\|e_T f\|_{W^1_p(\R_+;L_q(\Omega))} \leq 
2 \big( \|f\|_{L_p(0,T;W^2_q(\Omega))}+
\|f\|_{W^1_p(0,T;L_q(\Omega))}\big).
\end{equation}
Now we construct an extension $f^0$ of the initial data such that 
\begin{equation} \label{ext:init}
\|f^0\|_{L^p(0,T;W^2_q(\Omega))}
+\sup_{t \in \R_+}\|f^0\|_{B^{2-2/p}_{q,p}(\Omega)}
\leq C \|f(0)\|_{B^{2-2/p}_{q,p}(\Omega)}.
\end{equation}
and 
$$
f^0|_{\Omega \times \{t=0\}}=f(0,\cdot).
$$
For this purpose we can for instance extend $f(0,\cdot)$ to $g^0 \in B^{2-2/p}_{q,p}(W)$ such that $\Omega \subset W$ and 
\begin{equation} \label{gzero}
g^0|_{\de W}=0, \qquad \|g^0\|_{B^{2-2/p}_{q,p}(W)}\leq C(\Omega,V)\|f(0,\cdot)\|_{B^{2-2/p}_{q,p}(\Omega)}.
\end{equation}
Then we can find $f^0$ as a solution to the heat equation 
\begin{equation}
\begin{aligned}
&\de_t f^0 - \Delta f^0 = 0 \quad {\rm in} \quad W \times (0,T),\\
&f^0|_{\de V}=0, \quad f^0|_{t=0}=g^0.
\end{aligned}    
\end{equation}
Then the standard $L_p$-regularity theory of the heat equation yields 
$$
\|f^0\|_{L^p(0,T;W^2_q(W))}
+\sup_{t \in \R_+}\|f^0\|_{B^{2-2/p}_{q,p}(W)}
\leq C \|g^0\|_{B^{2-2/p}_{q,p}(W)},
$$
which together with \eqref{gzero} imply \eqref{ext:init}. 

Applying Lemma \ref{L:int} with $X=W^2_q(\Omega), \, Y=L_q(\Omega)$ and using 
\eqref{def:ext} and \eqref{ext:2} we have
\begin{align*}
&\sup_{t \in (0, T)}\|f(\cdot, t)-f^0\|_{B^{2(1-1/p)}_{q,p}(\Omega)}
\leq \sup_{t \in (0, \infty)}\|e_T[f-f^0]\|_{B^{2(1-1/p)}_{q,p}(\Omega)}\\&\quad 
\leq \Big(\|e_T[f-f^0]\|_{L_p((0, \infty), W^2_q(\Omega))}^p
+ \|e_T[f-f^0]\|_{W^1_p((0, \infty), L_q(\Omega))}^p\Big)^{1/p}\\
&\quad \leq C\Big(\|f-f^0\|_{L_p(0, \infty; W^2_q(\Omega))}
+ \|\pd_t f\|_{L_p(0, T; L_q(\Omega))}\Big).
\end{align*}
By \eqref{ext:init} we have 
\begin{align*}
\sup_{t \in (0,T)}\|f(\cdot, t)\|_{B^{2(1-1/p)}_{q,p}(\Omega)}&\leq
\sup_{t \in (0,T)}\|f(\cdot, t)-f^0\|_{B^{2(1-1/p)}_{q,p}(\Omega)} 
+\sup_{t \in (0,T)}\|f^0\|_{B^{2(1-1/p)}_{q,p}(\Omega)}\\[3pt]
&\leq C(\|f\|_{L_p(0, \infty; W^2_q(\Omega))}
+ \|\pd_t f\|_{L_p(0, T; L_q(\Omega))} + \|f(0)\|_{B^{2-2/p}_{q,p}(\Omega)}).
\end{align*}
This gives \eqref{4.1}, which together with Corollary \ref{c:imbed1} implies \eqref{4.2}.

\qed

\noindent
A direct conclusion from Lemma \ref{l:nonlin1} is 
\begin{lem} \label{l:nonlin2}
Assume $p,q$ satisfy the assumptions of Theorem \ref{t:local}. Assume moreover that
$$
\|(z_1,z_2,\vw)\|_{\cal X(T)}\leq M, \quad 
\|\vw(0)\|_{B^{2-2/p}_{q,p}(\Omega)} \leq L,\quad
\; z_1(0)=z_2(0)=0
$$
for some positive constants $M,L$.
Then 
\begin{align}
&\|\bV^0(\kv),
\nabla_{\kv}\bV^0(\kv)
\|_{L_\infty((0,T)\times\Omega)} \leq C(M,L)E(T), \label{est:01}\\
&{\rm sup}_{t \in (0,T)} \|z_i(\cdot,t)\|_{W^1_q(\Omega)}\leq C(M)E(T), \;i=1,2, \label{est:02} \\
&{\rm sup}_{t \in (0,T)}\|\vw(\cdot,t)-\vw(0)\|_{B^{2(1-1/p)}_{q,p}(\Omega)}\leq C(M,L),\label{est:03}\\
&\|\vw\|_{L_\infty(0,T;W^1_\infty(\Omega))}\leq C(M,L), \label{est:04}
\end{align}
where $\kv$ is defined in \eqref{kv}. 
\end{lem}
{\bf Proof.} \eqref{est:01} follows immediately from Lemma \ref{l:lag1}.
Next, we have 
$$\|z(\cdot, t)\|_{W^1_q(\Omega)}
\leq \int^t_0\|\de_t z(\cdot, s)\|_{W^1_q(\Omega)}\,ds
\leq T^{1/{p'}}\|\de_t z\|_{L_p((0, T), W^1_q(\Omega))}
\leq C(M)E(T),
$$
which implies \eqref{est:02}. Finally, \eqref{est:03} and \eqref{est:04} follow from \eqref{4.1} and \eqref{4.2}, respectively.  

\qed
\subsection{Estimate of the right hand side of \eqref{sysL2}}
Before going to the essence of this section let is observe the following fact
%
%
\begin{lem} \label{l:pos} Let $\fz=\fz(\fr,\fq)$ be as in \eqref{lag:5} and satisfy  \eqref{TZ} and \eqref{RleqZ}, then the following is true:
\begin{align}
&\fr+\fq>0 \; \Longrightarrow \; \fz>0, \label{pos1} \\
&\fr+\fq\geq \kappa>0 \; \Longrightarrow \; \fz \geq \min\left\{\frac{\kappa}{2},\left(\frac{\kappa}{2}\right)^{1/\gamma}\right\}, \label{pos2} \\ 
&\fz \leq \max\{ (2\fq)^{1/\gamma},2\fr \}. \label{pos3}
\end{align}
\end{lem}
{\em Proof.} If $\fq>0$ then $\fz>0$ by \eqref{TZ}, while if $\fr>0$ then $\fz>0$ by \eqref{RleqZ},
which proves \eqref{pos1}.

If $\fr+\fq\geq \kappa$ then $\fr\geq \frac{\kappa}{2}$ or $\fq\geq \frac{\kappa}{2}$. 
In the first case we have $\fz \geq \frac{\kappa}{2}$ by \eqref{RleqZ}. In the second case we have by \eqref{TZ} 
$$
\fz^\gamma=\left(1-\frac{\fr}{\fz}\right)^{-1}\fq \geq \fq \geq \frac{\kappa}{2},
$$
therefore we have \eqref{pos2}.

In order to show \eqref{pos3} we consider separately the cases $\fr \geq \frac{\fz}{2}$ and 
$\fr < \frac{\fz}{2}$. In the first case $\fz \geq 2\fr$. In the second we have $\left( 1-\frac{\fr}{\fz} \right)^{-1}<2$, which by \eqref{TZ} implies $\fz^\gamma<2\fq$. This completes the proof.

\qed

We are now ready to prove the following estimate for the RHS of \eqref{sysL2}: 
\begin{prop} \label{p:est:f} Assume $p,q$ satisfy the assumptions of Theorem \ref{t:local}.
Let $\bar V=(\bar \sigma,\bar \eta,\bar\vv)$, where $\bar\sigma=\bar\fr-\fr_0$, $\bar\eta=\bar\fq-\fq_0$ satisfy 
\begin{equation} \label{BM}
\|\bar V\|_{\CX(T)} \leq M, \quad \bar \vv|_{t=0}=\vv_0, \quad 
\bar\sigma+ \fr_0 + \bar \eta+ \fq_0\geq\delta,
\end{equation}
for some constants $M,\delta>0$. Then 
\begin{subequations} \label{est:f}
\begin{align}
&\|f_1(\bar V),f_2(\bar V)\|_{L_p(0,T;W^1_q(\Omega))} \leq C(M,L)E(T), \label{est:f1f2}\\
&\|\vf_3(\bar V)\|_{\lpqt} \leq C(M,L)E(T), \label{est:f3}
\end{align}
\end{subequations}
where $L$ is the constant from \eqref{init:loc}.
\end{prop}
{\bf Proof.} By \eqref{lag:6} and \eqref{f1} we have using summation convention: 
\begin{equation} \label{f1:1}
\begin{aligned}
f_1(\bar V)=-(\bar \sigma+ \fr_0)V_{ij}^0(\kbv)\frac{\de \bar v_i}{\de y_j}
-\bar \sigma \div \bar \vv,
\end{aligned}
\end{equation}
therefore by Lemma \ref{l:nonlin2} we immediately get 
\begin{equation}
\|f_1(\bar V)\|_{L_p(0,T;L_q(\Omega))}\leq C(M,L)E(T).    
\end{equation}
Differentiating \eqref{f1:1} we obtain
\begin{align*}
\nabla f_1(\bar V)=&-\nabla(\bar\sigma+ \fr_0)V_{ij}^0(\kbv)\frac{\de v_i}{\de y_j}
-(\bar\sigma+ \fr_0) \left[ \nabla_{\kv}V_{ij}^0(\kbv)\nabla\kbv\frac{\de v_i}{\de y_j}+V_{ij}^0(\kbv)\frac{\de^2 v_i}{\de y_i\de y_k} \right]\\
&-\div\bar\vv\nabla\bar\sigma-\bar\sigma\nabla\div \bar\vv.
\end{align*}
By Lemma \ref{l:nonlin2} we see that in each component one term is bounded in $L_p(0,T;L_q(\Omega))$ while all remaining terms are small in $L_\infty((0,T)\times \Omega)$, therefore 
\begin{equation}
\|\nabla f_1(\bar V)\|_{L_p(0,T;L_q(\Omega))}\leq C(M,L)E(T).    
\end{equation}
As $f_2(\bar V)$ has the same structure we obtain \eqref{est:f1f2}.
By \eqref{f3} we have
\begin{equation} \label{f3:2}
\begin{aligned}
\vf_3(\bar V)&=\vO_3(\bar V+(0, \fr_0, \fq_0))-(\bar \sigma+\bar \eta)\d_t\vv-I_1(\bar\fr,\bar\fq,\bar\fz)-I_2(\bar\fr,\bar\fq,\bar\fz)-I_3(\bar\fr,\bar\fq,\bar\fz)-I_4(\bar\fr,\bar\fq,\bar\fz), 
\end{aligned}
\end{equation}
where $\bar \fz$ is defined by $\bar\fq=\left( 1-\frac{\bar\fr}{\bar\fz} \right)\bar\fz^\gamma$.
Recalling \eqref{lag:7}, by Lemma \ref{l:nonlin2} we have 
\begin{equation} \label{5.0}
\|\vO_3(\bar V+(0, \fr_0, \fq_0)),(\bar \sigma+\bar \eta)\d_t\bar \vv\|_{L_p(0,T;L_q(\Omega))} \leq C(M,L)E(T).    
\end{equation}
In order to estimate the remaining parts observe first that all the denominators are bounded from below by positive powers of $\bar \fz$ and $\fz_0$. More precisely, 
\begin{equation} \label{5.1}
\begin{aligned}
&\gamma\fz_0-(\gamma-1) \fr_0 \geq \fz_0, & \gamma\bar\fz-(\gamma-1)\bar\fr \geq \bar\fz,\\
&\gamma\fz_0^{\gamma}-(\gamma-1) \fr_0\fz_0^{\gamma-1} \geq \fz_0^{\gamma},  &\gamma\bar\fz^{\gamma}-(\gamma-1)\bar\fr\,\bar\fz^{\gamma-1} \geq \bar\fz^{\gamma}.
\end{aligned}
\end{equation}
We have 
\begin{equation} \label{5.2}
\begin{aligned}
&\|\bar \fz(t) - \fz_0\|_{L_\infty(\Omega)}\leq \int_0^t \|\de_s\bar\fz(s)\|_{\linfsp}\,ds  
= \int_0^t \|\de_{\bar \sigma}\bar \fz \de_s \bar\sigma+\de_{\bar \fq}\bar \fz \de_s \bar\fq\|_{L_\infty(\Omega)}\,ds\\[3pt]
&\leq \|\de_{\bar\sigma} \bar\fz,\de_{\bar\fq} \bar\fz\|_{L_\infty(\Omega \times (0,T))} \int_0^t\|\de_s\bar\sigma,\de_t\bar\fq\|_{W^1_q(\Omega)}\ds\\[3pt]
&\leq C(\|\bar\fz,\fz_0,\bar\fz^{-1},\fz_0^{-1}\|_{{\linfsp}})t^{1/p'}\|\bar V\|_{\CX(T)}
\leq C(M,\|\bar\fz,\fz_0,\bar\fz^{-1},\fz_0^{-1}\|_{{\linfsp}})E(T).
\end{aligned}
\end{equation}
Similarly we get for any $\beta \in \R$ 
\begin{equation} \label{5.3}
\|\bar \fz^\beta(t) - \fz_0^\beta\|_{\linfsp} \leq
\int_0^t\|\bar\fz^{\beta-1}\de_s\bar\fz(s)\|_{\linfsp}\,ds 
\leq C(M,\|\bar\fz,\fz_0,\bar\fz^{-1},\fz_0^{-1}\|_{\linfsp})E(T).
\end{equation}
%
Now we can estimate $I_1$-$I_4$ defined in \eqref{f3}. By \eqref{est:02},\eqref{5.1} and \eqref{5.2} we have
\begin{equation} \label{5.4}
\begin{aligned}
&\|I_1(\bar \fr,\bar\fq,\bar\fz)(t)\|_{L_q(\Omega)} \leq C\|\bar\fz^{\gamma^+}-\fz_0^{\gamma^+},\bar\fz-\fz_0,\bar\fr- \fr_0\|_{L_\infty(\Omega\times(0,T))}\|\nabla\bar\sigma\|_{L_p(0,T;L_q)}\\[3pt]
&\leq C\big(M,\|\bar\fz,\fz_0,\bar\fz^{-1},\fz_0^{-1}\|_{\linfsp}\big)E(T)\|\bar V\|_{\CX(T)}\leq C\big(M,\|\bar\fz,\fz_0,\bar\fz^{-1},\fz_0^{-1}\|_{\linfsp}\big)E(T).
\end{aligned}
\end{equation}
%
Similarly, using also \eqref{5.3} with $\beta\in\{\gamma,\gamma-1,\gamma^+\}$ we obtain
\eq{ \label{I2I3}
&\|I_2(\bar \fr,\bar\fq,\bar\fz),I_3(\bar \fr,\bar\fq,\bar\fz),I_4(\bar \fr,\bar\fq,\bar\fz)\|_{\lpqt} \\
&\leq C\left(M, \|\fz_0^{-1},\bar\fz^{-1},\fz_0,\bar\fz\|_{\linf} \right)E(T),}
where smallness in time in the estimate for $I_4$ results from the fact that $\nabla  \fq_0$ and $\nabla \fr_0$ are time independent.
Now it is enough to observe that by Lemma \ref{l:pos} and \eqref{BM} we have
\begin{equation} \label{est:zeta}
\begin{aligned}
&\|\bar\fz\|_{\linf} \leq C (\|\bar\sigma+ \fr_0,\bar\eta+ \fq_0\|_{\linf} )= C(M,L),\\
&\|\fz_0\|_{\linf} \leq C (\| \fr_0, \fq_0\|_{\linf}) = C(M,L),\\
&\|\fz_0^{-1},\bar\fz^{-1}\|_{\linf} \leq C(M,L).
\end{aligned}
\end{equation}
Combining these bounds with \eqref{5.0}, \eqref{5.4} and \eqref{I2I3} we obtain \eqref{est:f3}.

\qed

\subsection{Contraction argument - proof of Theorem \ref{t:local}}

Let us define a solution operator\\ 
$$
(\sigma,\eta,\vv)=S(\bar \sigma, \bar \eta, \bar \vv)
$$
$$\iff$$
\centerline{$(\sigma,\eta,\vv)$ solves \eqref{sysL2} with the right hand side $(f_1(\bar V),f_2(\bar V),\vf_3(\bar V))$.}

\medskip
By \eqref{zeta:1}, as $\fz_0$ is bounded from below by a positive constant, we have 
$$
\|\nabla \fz_0\|_{L_q(\Omega)}\leq C\|\nabla \fr_0,\nabla \fq_0\|_{L_q(\Omega)},
$$
so by \eqref{gamma12:1} and \eqref{5.1} we see that 
$$
\|\nabla \omega_1^0,\nabla \omega_1^0\|_{L_q(\Omega)}\leq 
C\|\nabla \fz_0\|_{L_q(\Omega)}\leq C\|\nabla \fr_0,\nabla \fq_0\|_{L_q(\Omega)}.
$$
Therefore, by Theorem \ref{t:maxreg} and Proposition \ref{p:est:f}, $S$ is well defined on $\CX(T)$. 
First we show that for sufficiently small times $S$ maps certain bounded set in $\CX(T)$ into itself. Due to definition of $S(\cdot)$ we can restrict ourselves to functions satisfying the initial conditions of \eqref{sysL2}, this will be necessary for proving contractivity as we want to take advantage of smallness of time. Therefore we define 
\begin{equation} \label{def:BM}
B_M=\{(\bar\sigma,\bar\eta,\bar \vv) \in \CX(T):\; \|(\bar\sigma,\bar\eta,\bar \vv)\|_{\CX(T)}\leq M, \; (\bar\sigma,\bar\eta,\bar \vv)|_{t=0}=(0,0,\vv_0)\}.     
\end{equation}
\begin{lem}
For any $L>0$ there exists sufficiently large $M>0$ and sufficiently small $T>0$ such that 
$S(B_M) \subset B_M$, where $B_M$ is defined in \eqref{def:BM}.
\end{lem}
{\bf Proof:} By \eqref{est:linME1_A} and \eqref{est:f} we have 
$$
\|(\bar\sigma,\bar\eta,\bar\vv)\|_{\CX(T)} \leq M \; \Longrightarrow \;  \|S(\bar\sigma,\bar\eta,\bar\vv)\|_{\CX(T)} \leq C[L+C(M,L)E(T)],
$$
which gives the assertion. 

\qed
\noindent
Now we show that $S$ is a contraction on $B_M$. 
For this purpose set 
$$
V_i=(\sigma_i,\eta_i,\vv_i)=S(\bar V_i), \quad \bar V_i=(\bar \sigma_i,\bar \eta_i,\bar \vv_i), \quad i=1,2,  
$$
and define 
$$
\bar \fr_i = \fr_0+\bar \sigma_i, \quad \bar \fq_i = \fq_0+\bar \eta_i, \quad \bar \fq_i = \left( 1- \frac{\bar\fr_i}{\bar\fz_i}\right) \bar\fz_i^\gamma, \quad i=1,2. 
$$
Then the difference $(\delta \sigma,\delta \eta,\delta \vv)=V_1-V_2$ satisfies  
\eq{\label{sys:dif}
&\pt\delta\sigma+ \fr_0\Div\delta\vv=f_1(\bar V_1)-f_1(\bar V_2),\\
&\pt\delta\eta+ \fq_0\Div\delta\vv=f_2(\bar V_1)-f_2(\bar  V_2),\\
&( \fr_0+ \fq_0)\pt\delta\vv-\mu\lap\delta\vv-\nu\Grad\Div\delta\vv+\omega_1^0\Grad\delta\sigma+\omega_2^0\Grad\delta\eta=\vf_3(\bar V_1)-\vf_3(\bar V_2),\\
&(\delta\sigma,\delta\eta,\delta\vv)|_{t=0}=0, \quad \delta\vv|_{\d\Omega}=0.
}
By \eqref{est:linME1_A}, in order to prove that $S$ is a contraction on $B_M$ it is enough to show 
\begin{lem} \label{l:dif}
Assume $\bar V_1,\bar V_2 \subset B_M$. Then 
\begin{subequations}
\begin{align}
\|f_1(\bar V_1)-f_1(\bar V_2),f_2(\bar V_1)-f_2(\bar V_2)\|_{L_p(0,T;W^1_q(\Omega))} \leq C(M,L)E(T)\|\bar V_1-\bar V_2\|_{\CX(T)} \label{est:dif1}\\
\|\vf_3(\bar V_1)-\vf_3(\bar V_2)\|_{\lpqt} \leq C(M,L)E(T)\|\bar V_1-\bar V_2\|_{\CX(T)}. \label{est:dif2}
\end{align}
\end{subequations}
\end{lem}
{\bf Proof.} \eqref{est:dif1} can be shown analogously to \eqref{est:f1f2}. I order to show \eqref{est:dif2} we can follow the proof of \eqref{est:f3}. We directly get  
\begin{equation*}
\|\vO_3(\bar V_1)-\vO_3(\bar V_2),(\bar\sigma_1+\bar\eta_1)\de_t\bar\vv_1-(\bar\sigma_2+\bar\eta_2)\de_t\bar\vv_2\|_{\lpqt}\leq C(M,L)E(T)\|\bar V_1-\bar V_2\|_{\CX(T)} .   
\end{equation*}
With the remaining terms we shall go into some details due to more complicated structure. Let us focus on the estimate for $I_1(\bar \fr_1,\bar \fq_1,\bar\fz_1)-I_1(\bar \fr_2,\bar \fq_2,\bar\fz_2)$ (recall the definition \eqref{f3}). Let us denote 
$$
H_0=\frac{1}{\gamma\bar\fz_0-(\gamma-1) \fr_0}, \quad H_i = \frac{1}{\gamma\bar\fz_i-(\gamma-1)\bar\fr_i},\quad i=1,2.
$$
Then we have
\begin{equation*}
I_1(\bar \fr_1,\bar \fq_1,\bar\fz_1)-I_1(\bar \fr_2,\bar \fq_2,\bar\fz_2)
=\delta I_1^1+\delta I_1^2+\delta I_1^3+\delta I_1^4,
\end{equation*}
where 
\begin{equation} \label{deltaI}
\begin{aligned}
&\delta I_1^1=
\gamma\fz_0H_0 \,\Big[ H_1(\bar\fz_1^{\gamma^+}-\fz_0^{\gamma^+})\nabla(\bar\sigma_1-\bar\sigma_2)+(\bar\fz_1^{\gamma^+}-\fz_0^{\gamma^+})(H_1-H_2)\nabla\bar\sigma_2+(\bar\fz_1^{\gamma^+}-\bar\fz_2^{\gamma^+})H_2\nabla\bar\sigma_2 \Big],\\[3pt]
&\delta I_1^2=
\gamma\fz_0^{\gamma+}H_0 \,\Big[ H_1(\fz_0 -\bar\fz_1)\nabla(\bar\sigma_1-\bar\sigma_2)+(\fz_0-\bar\fz_1)(H_1-H_2)\nabla\bar\sigma_2+(\bar\fz_2-\bar\fz_1)H_2\nabla\bar\sigma_2 \Big]
\\[3pt]
&\delta I_1^3=
(\gamma-1)\fz_0^{\gamma+}H_0 \,\Big[ H_1(\bar\fr_1-\fr_0)\nabla(\bar\sigma_1-\bar\sigma_2)+(\bar\fr_1-\fr_0)(H_1-H_2)\nabla\bar\sigma_2+(\bar\fr_1-\bar\fr_2)H_2\nabla\bar\sigma_2 \Big],
\\[3pt]
&\delta I_1^4=
(\gamma-1)\fr_0H_0 \,\Big[ H_1(\fz_0^{\gamma^+}-\bar\fz_1^{\gamma^+})\nabla(\bar\sigma_1-\bar\sigma_2)+(\fz_0^{\gamma^+}-\bar\fz_1^{\gamma^+})(H_1-H_2)\nabla\bar\sigma_2+(\bar\fz_2^{\gamma^+}-\bar\fz_1^{\gamma^+})H_2\nabla\bar\sigma_2 \Big].
\end{aligned}
\end{equation}
Let us estimate $\delta I_1^1$, it has three components. First observe that estimates from the previous section imply 
\begin{equation} \label{deltaI:0}
\|\fz_0,H_0,H_1,H_2\|_{L_\infty(\Omega\times (0,T))} \leq C(M,L).    
\end{equation}
By \eqref{5.3} we have 
$$
\|\bar\fz_1^{\gamma^+}-\fz_0^{\gamma^+}\|_{L_\infty(\Omega \times (0,T))}\leq C(M)E(T),
$$
therefore 
\begin{equation} \label{deltaI:1}
\begin{aligned}
&\|(\bar\fz_1^{\gamma^+}-\fz_0^{\gamma^+})\nabla(\bar\sigma_1-\bar\sigma_2)\|_{\lpqt}
\leq \|(\bar\fz_1^{\gamma^+}-\fz_0^{\gamma^+})\|_{L_\infty(\Omega \times (0,T))}
\|\nabla(\bar\sigma_1-\bar\sigma_2)\|_{\lpqt}\\[3pt]
&\leq C(M)E(T)\|\bar V_1-\bar V_2\|_{\CX(T)}.
\end{aligned}
\end{equation}
For the second term we have 
\begin{equation} \label{deltaI2:a}
\begin{aligned}
&\|(\bar\fz_1^{\gamma^+}-\fz_0^{\gamma^+})(H_1-H_2)\nabla\bar\sigma_2\|_{\lpqt}\\[3pt]
&\leq \|(\bar\fz_1^{\gamma^+}-\fz_0^{\gamma^+})\|_{L_\infty(\Omega\times(0,T))}\|\nabla\bar\sigma_2\|_{\lpqt}\|H_1-H_2\|_{L_\infty(\Omega\times (0,T))}\\[3pt]
&\leq C(M)E(T)\|H_1-H_2\|_{L_\infty(\Omega\times (0,T))}.
\end{aligned}
\end{equation}
We have
$$
H_1-H_2=\frac{ \gamma(\bzeta_2-\bzeta_1)-(\gamma-1)(\bar\fr_2-\bar\fr_1) }{ [\gamma\bzeta_1-(\gamma-1)\bar\fr_1][\gamma\bzeta_2-(\gamma-1)\bar\fr_2]}.
$$
Therefore, by \eqref{def:BM} and the estimates from the previous section,
$$
\|H_1-H_2\|_{L_\infty(\Omega\times(0,T))} \leq C(M) \|\bzeta_1-\bzeta_2,\bvr_1-\bvr_2\|_{L_\infty(\Omega\times(0,T))}
\leq C(M)E(T)\|\bar V_1-\bar V_2\|_{\CX(T)},
$$
which together with \eqref{deltaI2:a} implies
\begin{equation} \label{deltaI:2}
\|(\bar\fz_1^{\gamma^+}-\fz_0^{\gamma^+})(H_1-H_2)\nabla\bar\sigma_2\|_{\lpqt}
\leq C(M)E(T)\|\bar V_1-\bar V_2\|_{\CX(T)}.
\end{equation}
And finally the third term. By \eqref{def:BM} we have $(\bar\fz_1-\bar\fz_2)|_{t=0}=0$, therefore 
\begin{equation*}
\begin{aligned}
&\|(\bar\fz_1^{\gamma^+}-\bar\fz_2^{\gamma^+})H_2\nabla\bar\sigma_2\|_{\lpqt}
\leq C\|\bar\fz_1^{\gamma^+}-\bar\fz_2^{\gamma^+}\|_{L_\infty(\Omega\times (0,T))}\|\nabla\bar\sigma_2\|_{\lpqt}\\[3pt]
&\leq C(M) \int_0^t \|\de_t(\bar\fz_1^{\gamma^+}-\bar\fz_2^{\gamma^+})(s)\|_{\linfsp}\,ds\\[3pt]
&\leq C(M) \int_0^t \Big[ \|(\bar\fz_1^{\gamma^+-1}-\bar\fz_2^{\gamma^+-1})\de_t\bar\fz_2\|_{\linfsp}+\|\bar\fz_1^{\gamma^+-1}\de_t(\bar\fz_1-\bar\fz_2)\|_{\linfsp}\Big]\,ds\\[3pt]
&\leq C(M)\|\bar\fz_1^{\gamma^+-1}-\bar\fz_2^{\gamma^+-1}\|_{L_\infty(\Omega\times(0,T))}\int_0^t\|\de_t\bar\fz_2\|_{\linfsp}\,ds 
+\|\bar\fz_1^{\gamma^+-1}\|_{L_\infty(\Omega\times(0,T))}\int_0^t\|\de_t(\bar\fz_1-\bar\fz_2)\|_{\linfsp}\,ds\\[3pt]
&\leq C(M)\|\bar\sigma_1-\bar\sigma_2,\bar\eta_1-\bar\eta_2\|_{L_\infty(\Omega\times(0,T))}E(T)+C(M)E(T)\|\bar V_1-\bar V_2\|_{\CX(T)}.
\end{aligned}
\end{equation*}
But $(\bar\sigma_1-\bar\sigma_2,\bar\eta_1-\bar\eta_2)|_{t=0}=(0,0)$, therefore 
$$
\|\sigma_1-\sigma_2,\eta_1-\eta_2\|_{L_\infty(\Omega\times(0,T))}\leq E(T)\|\bar V_1-\bar V_2\|_{\CX(T)},
$$
so altogether 
\begin{equation} \label{deltaI:3}
\|(\fz_1^{\gamma^+}-\fz_2^{\gamma^+})H_2\nabla\sigma_2\|_{\lpqt} \leq C(M)E(T)\|\bar V_1-\bar V_2\|_{\CX(T)}.    
\end{equation}
Combining \eqref{deltaI:0},\eqref{deltaI:1},\eqref{deltaI:2} and \eqref{deltaI:3} we obtain 
\begin{equation} \label{deltaI1}
\|\delta I_1^1\|_{\lpqt}  \leq C(M)E(T)\|\bar V_1-\bar V_2\|_{\CX(T)}.
\end{equation}
The terms $\delta I_2-\delta I_4$ have similar structure to $\delta I_1$ therefore we estimate them in the same way obtaining altogether   
\begin{equation*}
\|\delta I_1\|_{\lpqt}  \leq C(M)E(T)\|\bar V_1-\bar V_2\|_{\CX(T)}.
\end{equation*}
Now, if we would like to write precise form of the terms 
$I_2(\bar \fr_1,\bar \fq_1,\bar\fz_1)-I_2(\bar \fr_2,\bar \fq_2,\bar\fz_2)$, 
$I_3(\bar \fr_1,\bar \fq_1,\bar\fz_1)-I_3(\bar \fr_2,\bar \fq_2,\bar\fz_2)$
and
$I_4(\bar \fr_1,\bar \fq_1,\bar\fz_1)-I_4(\bar \fr_2,\bar \fq_2,\bar\fz_2)$
it would be convenient to define 
$$
H_{0,\gamma}=\frac{1}{\gamma\bar\fz_0^\gamma-(\gamma-1) \fr_0\bar\fz_0^{\gamma-1}}, \quad  H_{i,\gamma} = \frac{1}{\gamma\bar\fz_i^\gamma-(\gamma-1)\bar\fr_i\bar\fz_i^{\gamma-1}},\quad i=1,2.
$$
Then we obtain expressions with structure similar to \eqref{deltaI} with functions $ H_{0,\gamma}, H_{1,\gamma}, H_{2,\gamma}$ instead of $H_0, H_1, H_2$. Now it is enough to observe that estimates from the previous section imply 

\begin{equation} \label{deltaI:0b}
\|H_{0,\gamma}, H_{1,\gamma}, H_{2,\gamma}\|_{L_\infty(\Omega\times (0,T))} \leq C(M,L)    
\end{equation}
and 
$$
\|H_{1,\gamma}-H_{2,\gamma}\|_{L_\infty(\Omega\times(0,T))} \leq
C(M)E(T)\|\bar V_1-\bar V_2\|_{\CX(T)},
$$
since 
$$
H_{1,\gamma}- H_{2,\gamma} = \frac{\gamma(\bar \fz_1^\gamma-\bar \fz_2^\gamma)-(\gamma-1)[ (\bar r_2-\bar r_1)\bar\fz^{\gamma-1} +\bar r_1(\bar\fz_2^{\gamma-1}-\bar\fz_1^{\gamma-1}) ] }{[\gamma\bar\fz_1^\gamma-(\gamma-1) \fr_1\bar\fz_1^{\gamma-1}][\gamma\bar\fz_2^\gamma-(\gamma-1) \fr_2\bar\fz_2^{\gamma-1}]}.
$$
Therefore we prefer not to bother the reader with repeating estimates similar to those leading to \eqref{deltaI1}  and notice only that finally we obtain  
\begin{equation*}
\begin{aligned}
&\|I_2(\bar \fr_1,\bar\fq_1,\bar\fz_1)-I_1(\bar \fr_2,\bar\fq_2,\bar\fz_2),
I_3(\bar \fr_1,\bar\fq_1,\bar\fz_1)-I_3(\bar \fr_2,\bar\fq_2,\bar\fz_2),I_4(\bar \fr_1,\bar\fq_1,\bar\fz_1)-I_4(\bar \fr_2,\bar\fq_2,\bar\fz_2)\|_{\lpqt}\\ 
&\leq C(M,L)E(T)\|\bar V_1-\bar V_2\|_{\CX(T)},
\end{aligned}
\end{equation*}
so altogether we have \eqref{est:dif2}.

\qed
\noindent
Now, applying Theorem \ref{t:maxreg} and Lemma \ref{l:dif} to \eqref{sys:dif} we obtain 
$$
\|V_1-V_2\|_{\CX(T)} \leq C(M,L)E(T) \|\bar V_1-\bar V_2\|_{\CX(T)},
$$
and Theorem \ref{t:local} follows from the Banach contraction principle.

\section{Global well-posedness} \label{s:glob}
\subsection{Linearization around the constant state}
In order to prove Theorem \ref{t:global} we linearize \eqref{sysL1} around the positive constants $R_*,Q_*$ to which the initial condition is assumed to be close. For convenience let us denote 
$$
\fq_*=Q_*, \quad \fr_*=R_*.
$$
Then we define  
\begin{equation*}
\sigma=\fr-\fr_*, \quad \eta=\fq-\fq_*    
\end{equation*}
and 
\begin{equation*}
V=(\vv,\sigma,\eta)=U-(0,\fr_*,\fq_*), 
\end{equation*}
where $U$ is defined in \eqref{def:U} (notice that we use the same notation for perturbations as in the previous section, the reason is we prefer to avoid introducing additional notation). We obtain 
\eq{\label{sysL2:const}
&\pt\sigma+\fr_*\Div\vv=g_1(V),\\[3pt]
&\pt\eta+\fq_*\Div\vv=g_2(V),\\[3pt]
&(\fr_*+\fq_*)\pt\vv-\mu\lap\vv-\nu\Grad\Div\vv+\omega_{1*}\Grad\sigma+\omega_{2*}\Grad\eta=\vg_3(V),\\[3pt]
&(\sigma,\eta,\vv)|_{t=0}=( \fr_0-\fr_*, \fq_0-\fq_*,\vv_0), \quad \vv|_{\d\Omega}=0,
}
where $\fz_*=\fz_*(\fr_*,\fq_*)$ is a constant defined by 
\begin{equation*} 
\fq_*=\left(1-\frac{\fr_*}{\fz_*} \right)\fz_*^\gamma    
\end{equation*}
and, analogously to \eqref{gamma12:1},
\begin{equation} \label{gamma12:2} 
\omega_{1*} = \frac{\gamma^+\fz_*^{\gamma^+}}{\gamma\fz_*-(\gamma-1)\fr_*}, \qquad
\omega_{2*} = \frac{\gamma^+\fz_*^{\gamma^+}}{\gamma\fz_*^{\gamma}-(\gamma-1)\fr_*\fz_*^{\gamma-1}}.
\end{equation}
The right hand side of \eqref{sysL2:const} is analogous to \eqref{f1}-\eqref{f3} with $(\fr_*,\fq_*,\fz_*)$ instead of $( \fr_0, \fq_0,\fz_0)$ and without terms with $\nabla  \fr_0$ and $\nabla  \fq_0$ in $\vg_3$: 
\begin{align}
g_1(V)&=O_1(V+(0,\fr_*,\fq_*))-\sigma\Div\vv, \label{f1a} \\[5pt]
g_2(V)&=O_2(V+(0,\fr_*,\fq_*))-\eta\Div\vv, \label{f2a} \\[8pt]
\vg_3(V)&=\vO_3(V+(0,\fr_*,\fq_*))-(\sigma+\eta)\d_t\vv \label{f3a} \\[5pt]
\nonumber&-\underbrace{ \frac{ \gamma^+\gamma\big[\fz_*( \fz^{\gamma^+}-\fz_*^{\gamma^+} )+\fz_*^{\gamma^+}(\fz_*-\fz)\big]+(\gamma-1)\big[ \fz_*^{\gamma^+}(\fr-\fr_*)+\fr_*(\fz_*^{\gamma^+}-\fz^{\gamma^+}) \big] }{\big(\gamma\fz_*-(\gamma-1)\fr_*\big)\big(\gamma\fz-(\gamma-1)\fr\big) }\nabla \sigma }_{J_1(\fr,\fq,\fz)} \\[5pt]
\nonumber& - \underbrace{ \frac{ \gamma^+\gamma\big[ \fz_*^{\gamma^+}(\fz_*^{\gamma}-\fz^{\gamma})+\fz_*^{\gamma}(\fz^{\gamma^+}-\fz_*^{\gamma^+}) \big] }{ \big(\gamma\fz^{\gamma}-(\gamma-1)\fr\fz^{\gamma-1}\big)\big(\gamma\fz_*^{\gamma}-(\gamma-1)\fr_*\fz_*^{\gamma-1}\big)}\nabla\eta}_{J_2(\fr,\fq,\fz)}\\[5pt]
\nonumber& -\underbrace{\frac{\gamma^+(\gamma-1)\big[ \fr\fz_*^{\gamma^+}(\fz^{\gamma-1}-\fz_*^{\gamma-1})+\fz_*^{\gamma^+}\fz_*^{\gamma-1}(\fr-\fr_*)+\fr_*\fz_*^{\gamma-1}(\fz_*^{\gamma^+}-\fz^{\gamma^+}) \big]  }{ \big(\gamma\fz^{\gamma}-(\gamma-1)\fr\fz^{\gamma-1}\big)\big(\gamma\fz_*^{\gamma}-(\gamma-1)\fr_*\fz_*^{\gamma-1}\big)}\nabla\eta}_{J_3(\fr,\fq,\fz)}.  
\end{align}

\subsection{Exponential decay} \label{s:decay}
Consider a linear problem corresponding to \eqref{sysL2:const}
\eq{\label{sys:lin2}
&\pt\sigma+\fr_*\Div\vv=f_1,\\[3pt]
&\pt\eta+\fq_*\Div\vv=f_2,\\[3pt]
&(\fr_*+\fq_*)\pt\vv-\mu\lap\vv-\nu\Grad\Div\vv+\omega_{1*}\Grad\sigma+\omega_{2*}\Grad\eta=\vf_3,\\[3pt]
&(\sigma,\eta,\vv)|_{t=0}=(\sigma_0,\eta_0,\vv_0), \quad \vv|_{\d\Omega}=0.
}
For the above problem we have exponential decay result  
\begin{thm} \label{t:decay}
Let $p,q,\sigma_0,\eta_0$ and $\vv_0$ satisfy the assumptions of
Theorem \ref{t:maxreg}. Assume moreover that 
\begin{equation*}
e^{\beta t}f_1,e^{\beta t}f_2 \in L_p(\R_+,W^1_q(\Omega)), e^{\beta t}\vc{f}_3 \in L_p(\R_+,L_q(\Omega))
\end{equation*}
for some constant $\beta>0$.
Then solution to \eqref{sys:lin2} exists globally in time and satisfies 
\eq{ \label{est:decay}
\|e^{\beta t}\sigma,e^{\beta t}\eta,e^{\beta t}\vv\|_{\dot{\cal X}(+\infty)} \leq & C_{p,q}\, [\|\vv_0\|_{B^{2-2/p}_{q,p}(\Omega)}
+\|\sigma_0,\eta_0\|_{W^1_q(\Omega)}\\
&+\|e^{\beta t}\vf_3\|_{L_p(0,T;L_q(\Omega)}+\|(e^{\beta t}f_1,e^{\beta t}f_2)\|_{L_p(0,T;W^1_q(\Omega)}].
}
\end{thm}
The keynote ideas necessary to prove Theorem \ref{t:decay} were developed by Enomoto and Shibata \cite[Theorem 2.9]{ES}, where analogous result for the compressible Stokes system 
was shown, however with additional zero mean condition for the right hand side of the continuity equation ($f_1,f_2$ in our case). This condition was removed in \cite[Theorem 5.1]{PSZ}, where a system describing a mixture of two constituents is considered and its linearization consists of two parabolic equations coupled with the continuity equation. It is therefore slightly more complicated then the compressible Stokes system, but the essential ideas from \cite{ES} apply. In fact, the zero mean condition is used at the first stage of the proof \cite[Theorem 7.1]{PSZ} to guarantee the unique solvability of the resolvent problem for $\lambda=0$, but then the general case is reduced to the zero mean case. The details can be found in Section 7 of \cite{PSZ}.

Analogous result for the compressible Stokes system was proved recently in \cite[Proposition 3.2]{KNP}, also without the zero mean assumption. The proof is based on the ideas from \cite{PSZ} but the system is simpler, so we recommend also this proof for more details. 

Theorem \eqref{t:decay} can be proved following line by line the proof of Proposition 3.2 from \cite{KNP}, therefore we omit the proof here.

\subsection{Bounds for nonlinearities}
We start with the following counterpart of Lemma \ref{l:nonlin2}:
\begin{lem} \label{l:nonlin3}
Let $(\ebt z_1, \ebt z_2,\ebt \vw) \in \dot{\cal X}(\infty)$ and 
$$
\|z_1(0)-\fr_*,z_2(0)-\fq_*\|_{W^1_q(\Omega)}  +\|\vw(0)\|_{B^{2-2/p}_{q,p}(\Omega)} \leq \ep.
$$
Then 
\begin{align}
&\|\bV^0(\kw),
\nabla_{\kw}\bV^0(\kw)
\|_{L_\infty((0,T)\times\Omega)} \leq C\|\ebt (z_1,z_2,\vw) \|_{\dot \CX(+\infty)}, \label{est:01a}\\
&{\rm sup}_{t \in \R_+} \|z_1(\cdot,t)-\fr_*\|_{W^1_q(\Omega)}\leq C(\ep+\|\ebt (z_1,z_2,\vw) \|_{\dot \CX(+\infty)}),  \label{est:02a} \\
&{\rm sup}_{t \in \R_+} \|z_2(\cdot,t)-\fq_*\|_{W^1_q(\Omega)}\leq C(\ep+\|\ebt (z_1,z_2,\vw) \|_{\dot \CX(+\infty)}),  \label{est:02b} \\
&{\rm sup}_{t \in \R_+}\|\vw(\cdot,t)-\vw(0)\|_{B^{2(1-1/p)}_{q,p}(\Omega)}\leq C\|\ebt (z_1,z_2,\vw) \|_{\dot \CX(+\infty)},\label{est:03a}\\
&\|\vw\|_{L_\infty(0,T;W^1_\infty(\Omega))}\leq C\|\ebt (z_1,z_2,\vw) \|_{\dot \CX(+\infty)}, \label{est:04a}
\end{align}
where $\kw$ is defined as in \eqref{kv}.
\end{lem}
{\bf Proof.} We have
\begin{equation}
\begin{aligned}
\|\bV^0(\kw),
&\nabla_{\kw}\bV^0(\kw)\|_{L_\infty((0,T)\times\Omega)} \leq
C(\bV^0)\int_0^\infty \|\nabla_y \vw\|_{\infty} \dt\\
&\leq \left(\int_0^\infty e^{-\beta t p'} \dt\right)^{1/p'} \left(\int_0^\infty e^{\beta t p}\|\vw\|_{W^2_q(\Omega_0)} \dt\right)^{1/p},
\end{aligned}
\end{equation}

which implies \eqref{est:01a}. Next, 
\eqh{
&\|z_1(\cdot,t)-\fr_*\|_{W^1_q(\Omega_0)}\leq \|z_1(0)-\fr_*\|_{W^1_q(\Omega_0)}+\int_0^t\|\de_tz_1(s,\cdot)\|_{W^1_q(\Omega_0)} \dt\\
&\leq \ep+ C\left(\int_0^t e^{-\beta s p'} \,ds\right)^{1/p'} \left(\int_0^\infty e^{\beta s p}\|z_t\|_{W^1_q(\Omega_0)} \,d s\right)^{1/p},
}
and similarly for $z_2$ and $\fq_*$, which gives \eqref{est:02a} and \eqref{est:02b}. 
Finally, \eqref{est:03a} and \eqref{est:04a} follows from Lemma \ref{l:nonlin1}.

\qed

We are now in a position to prove the following estimate for the RHS of \eqref{sysL2}, which will enable us to apply Theorem \ref{t:decay}: 
\begin{prop} \label{p:nonlin2}
Let $g_1(V), g_2(V),\vc{g}_3(V)$ be defined by \eqref{f1a}-\eqref{f3a}.
Then
\begin{equation} \label{6.0}
\begin{aligned}
&\|\ebt g_1(V),\ebt g_2(V)\|_{L_p(\R_+,W^1_q(\Omega))}
+\|\ebt \vg_3(V)\|_{L_p(\R_+,L_q(\Omega))} \\[3pt] 
&\leq C \Big( \| \ebt V \|_{\dot \CX(+\infty)}^2 + \ep\big(\| \ebt V \|_{\dot \CX(+\infty)}+1\big)  \Big),
\end{aligned}
\end{equation}
where $\ep$ is from \eqref{small}.
\end{prop}
{\bf Proof.} By \eqref{lag:6} and \eqref{f1a} we have 
\begin{equation} \label{f1:1a}
\begin{aligned}
g_1(V)=-(\sigma+\fr_*)V_{ij}^0(\kv)\frac{\de v_i}{\de y_j}
- \sigma \div \bv,
\end{aligned}
\end{equation}
therefore by Lemma \ref{l:nonlin3} we immediately get 
\begin{equation*}
\|\ebt g_1(V)\|_{L_p(0,T;L_q(\Omega))}\leq C\|\ebt V\|^2_{\dot \CX(+\infty)}. \end{equation*}
Differentiating \eqref{f1:1a} we obtain
\begin{align*}
\nabla g_1(V)=&-\nabla(\sigma V_{ij}^0(\kv)\frac{\de v_i}{\de y_j}
-(\sigma+\fr_*) \left[ \nabla_{\kv}V_{ij}^0(\kv)\nabla\kv\frac{\de v_i}{\de y_j}+V_{ij}^0(\kv)\frac{\de^2 v_i}{\de y_i\de y_k} \right]\\
&-\div \vv\nabla\sigma-\sigma\nabla\div \vv.
\end{align*}
By Lemma \ref{l:nonlin3} we easily obtain  
\begin{equation*}
\|\ebt \nabla g_1(V)\|_{L_p(0,T;L_q(\Omega))}\leq C\|\ebt V\|^2_{\dot \CX(+\infty)}.
\end{equation*}
As $g_2$ has the same form with $ \eta,\fq_*$ instead of $\sigma,\fr_*$, alltogether we get 
\begin{equation} \label{6.1}
\|\ebt g_1(V),\ebt g_2(V)\|_{L_p(\R_+,W^1_q(\Omega))}
\leq 
C \| \ebt V \|_{\dot \CX(+\infty)}^2.    
\end{equation}
It remains to estimate $\vg_3$. This is more involved, however we can take advantage of ideas from the proof of local well posedness, roughly speaking replacing smallness of time by smallness of initial data. Firstly, by 
\eqref{lag:7} and Lemma \ref{l:nonlin3} we have 
\begin{equation} \label{6.2}
\|\ebt\vO_3( V+(0,\fr_*,\fq_*)),\ebt( \sigma+ \eta)\d_t \vv\|_{L_p(0,T;L_q(\Omega))} \leq C \| \ebt V \|_{\dot \CX(+\infty)}^2.     
\end{equation}
Now we have to estimate $J_1, J_2,J_3$. Taking into account \eqref{5.1} with $(\fr_*,\fq_*,\fz_*)$ instead of $( \fr_0,\fq_0,\fz_0)$, it is enough to find bounds on the numerators of the last three terms of \eqref{f3a}. 
These contains many terms but of a similar structure. First representative case is 
\begin{align*}
&\|\ebt (\fr-\fr_*)\nabla \eta\|_{L_p(\R_+;L_q(\Omega))} 
\leq C \|\fr-\fr_*\|_{L_\infty(\Omega\times\R_+)}\|\ebt \nabla\eta\|_{L_p(\R_+;L_q(\Omega))}\\[3pt] 
&\leq C \| \ebt V \|_{\dot \CX(+\infty)}\big( \| \ebt V \|_{\dot \CX(+\infty)} + \ep \big),
\end{align*}
where we have used \eqref{est:02a}.
The second case is when we have increment of $\fz$. Then we use
\begin{align*}
&\|\fz(t)-\fz_*\|_{\linfsp}\leq \|\fz_0-\fz_*\|_{\infty}+\|\fz-\fz_0\|_{\linfsp}\\[5pt]
&\leq \|( \fr_0-\fr_*)\de_{\vr}\fz,( \fq_0-\fq_*)\de_{\kappa}\fz\|_{\linfsp}+\int_0^t\|\de_t\fz(s)\|_{W^1_q}\,ds 
\leq C\big( \ep+\|\ebt V\|_{\dot\CX(t)} \big)
\end{align*}
to obtain 
\begin{align*}
&\|\ebt\fz^{\gamma^+}(\fz_*-\fz)\nabla\sigma\|_{L_p(\R_+;L_q(\Omega))}
\leq C\|\fz-\fz_*\|_{L_\infty(\Omega \times \R_+)}\|\ebt\nabla\sigma\|_{L_p(\R_+;L_q(\Omega))}\\[3pt] 
&\leq C\big( \ep+\|\ebt V\|_{\dot\CX(+\infty)} \big)\|\ebt V\|_{\dot\CX(+\infty)}.
\end{align*}
Other terms containing $(\fz_*-\fz)$ can be treated similarly. 
Terms with $\fz^\beta-\fz_*^\beta$ for \\$\beta \in \{\gamma^+,\gamma^+-1,\gamma-1 \}$ can also be estimated in a similar way. Altogether we obtain 
\begin{equation} \label{6.3}
\sum_{k=1}^3\| \ebt J_k(\fr,\fq,\fz)\|_{L_p(\R_+;L_q(\Omega))}
\leq C\big( \ep+\|\ebt V\|_{\dot\CX(+\infty)} \big)\|\ebt V\|_{\dot\CX(+\infty)}.
\end{equation}
Putting together \eqref{6.1},\eqref{6.2} and \eqref{6.3} we obtain \eqref{6.0}.

\qed

\subsection{Proof of Theorem \ref{t:global}}
We proceed in a standard way, first we show 
\begin{lem} Assume $(\sigma,\eta,\vv)$ is solution to \eqref{sysL2} with $ \fr_0, \fq_0$ and $\vu_0$ satisfying the assumptions of Theorem \ref{t:global}. Then 
\begin{equation} \label{est:expnorm}
\|\ebt(\sigma,\eta,\vv)\|_{\dot \CX(+\infty)} \leq E(\epsilon).    
\end{equation}
\end{lem}
{\bf Proof.} Applying Theorem \ref{t:decay} and Proposition \ref{p:nonlin2} to \eqref{sysL2} we obtain 
\begin{equation} \label{est:expnorm1}
\|\ebt(\sigma,\eta,\vv)\|_{\dot \CX(T)} \leq C \Big( \ep + \|\ebt(\sigma,\eta,\vv)\|_{\dot \CX(T)}^2 \Big) \quad \forall \; 0<T\leq +\infty.
\end{equation}
Notice that we derived this inequality for $T = \infty$, but the same arguments allow to obtain it for any $T > 0$. Now consider the equation 
$$
x^2-\frac{x}{C}+\epsilon =0 
$$
with roots  
$$
x_1(\epsilon)= \frac{1}{2C}-\sqrt{\frac{1}{4C^2}-\epsilon},
\quad x_2(\epsilon)= \frac{1}{2C}+\sqrt{\frac{1}{4C^2}-\epsilon}. 
$$
Observe that inequality \eqref{est:expnorm1} implies either $\|\ebt (\sigma,\eta,\vv)\|_{\dot \CX(T)} \leq x_1(\epsilon)$
or $\|\ebt (\sigma,\eta,\vv)\|_{\dot \CX(T)} \geq x_2(\epsilon)$. However, 
$$
\lim_{T \to 0}\|\ebt (\sigma,\eta,\vv)\|_{\dot \CX(T)} = 0,
$$
therefore 
\begin{equation} \label{6:3}
\|\ebt (\sigma,\eta,\vv)\|_{\dot \CX(T)} \leq x_1(\epsilon)
\end{equation}
for small times. However, $\|\ebt (\sigma,\eta,\vv)\|_{\dot \CX(T)}$ is continuous in time and therefore \eqref{6:3} holds for $0<T\leq +\infty$.

\qed

Now we can prolong the local solution for arbitrarily large times. For this purpose observe that if the initial data satisfies \eqref{small} then the time of existence from Theorem \ref{t:local} satisfies $T>C(\epsilon)>0$. Therefore, 
for arbitrarily large $T^*$ we obtain a solution on $(0,T^*)$ in a finite number of steps. By the estimate \eqref{est:expnorm} this solution satisfies \eqref{est:glob}.

\section*{Appendix}
Here we recall the definition of $\CR$-boundedness:
\begin{df}\label{dfn:7.1} Let $X$ and $Y$ be two Banach spaces, and 
$\|\cdot\|_X$ and $\|\cdot\|_Y$ their norms.  
A family of operators $\CT \subset \CL(X, Y)$ is called 
$\CR$-bounded on $\CL(X, Y)$ if there exist constants $C > 0$
and $p \in [1, \infty)$ such that 
for any $n \in \N$, $\{T_j\}_{j=1}^n \subset \CT$ and $\{f_j\}_{j=1}^n
\subset X$, the inequality 
$$\int^1_0\|\sum_{j=1}^n r_j(u)T_jf_j\|_Y^p\,du
\leq C\int^1_0\|\sum_{j=1}^nr_j(u)f_j\|_X^p\,du,
$$
where $r_j: [0, 1] \to \{-1,1\}$, $j \in \N$, are the Rademacher functions
given by 
$r_j(t) = {\rm sign}(\sin (2^j\pi t))$.  The smallest such
$C$ is called $\CR$-bound of $\CT$ on 
$\CL(X, Y)$ which is written by $\CR_{\CL(X, Y)}\CT$.
\end{df}
Next we quote Weis' vector valued Fourier Multiplier Theorem (\cite{Weis}):
\begin{thm} \label{thm:Weis}
Let $X$ and $Y$ be UMD spaces and $1<p<\infty$. Let $M \in C^1(\R\setminus \{0\}, \CL(X,Y))$. Let us define the operator 
$T_M:\CF^{-1} \CD(\R,X) \to \CS'(\R,Y)$:
\begin{equation} \label{def:TM}
T_M\phi(\tau) = \CF^{-1}[M\CF[\phi](\tau)].    
\end{equation}
Assume that 
\begin{equation} \label{rbound:Weis}
\CR_{\CL(X,Y)}(\{M(\tau): \; \tau \in \R \setminus \{0\}\})= \kappa_0<\infty, \quad 
\CR_{\CL(X,Y)}(\{\tau M'(\tau): \; \tau \in \R \setminus \{0\}\})=\kappa_1<\infty.
\end{equation}
Then, the operator $T_M$ defined in \eqref{def:TM} is extended to a bounded linear operator $L_p(\R,X) \rightarrow L_p(\R,X)$ and 
$$
\|T_M\|_{\CL(L_p(\R,X),L_p(\R,Y))} \leq C( \kappa_0+\kappa_1),
$$
where $C=C(p,X,Y)>0$. 
\end{thm}
\begin{rmk}
For definitions and properties of UMD spaces we refer the reader for example to Chapter 4 in \cite{HNVW}. Here let us only note that $L_p$ spaces and $W^k_p$ spaces are UMD for $1<p<\infty$.
\end{rmk}
We recall that the Fourier transform and its inverse are defined as  
\begin{equation} \label{def:FT}
\CF[f](\tau)=\int_\R e^{-it\tau}f(t)dt, \quad \CF^{-1}[f](t)=\frac{1}{2\pi}\int_\R e^{it\tau}f(\tau)d\tau   
\end{equation}
while the Laplace transform and its inverse are 
\begin{equation} \label{def:LT}
\CL[f](\lambda)=\int_\R e^{-\lambda t}f(t)dt, \quad \CL^{-1}_\beta[f](t)=\frac{1}{2\pi}\int_\R e^{\lambda t}f(\lambda)d\tau, \quad \textrm{where} \; \lambda=\beta+i\tau.   
\end{equation}
The following result (Theorem 2.17 in \cite{ES}) explains how to obtain $L_p$ - maximal regularity using Theorem \ref{thm:Weis} and Laplace transform:  
\begin{thm} \label{thm:maxreg}
Let $X$ and $Y$ be UMD Banach spaces and $1<p<\infty$. Let $0<\ep<\frac{\pi}{2}$ and $\beta_1 \in \R$. 
Let $\Phi_\lambda$ be a $C^1$ function of $\tau \in \R\setminus\{0\}$ where $\lambda=\beta+i\tau \in \Sigma_{\ep}+\beta_1$ with values in $\CL(X,Y)$. 
Assume that 
\begin{equation*}  
\CR_{\CL(X,Y)}(\{\Phi_\lambda: \; \lambda \in \Sigma_{\ep}+\beta_1\})\leq M,\quad
\CR_{\CL(X,Y)}\left(\left\{ \tau\frac{\de}{\de\tau}\Phi_\lambda: \; \lambda \in \Sigma_{\ep}+\beta_1\right\}\right) \leq M
\end{equation*}
for some $M>0$. Let us define 
\begin{equation} \label{def:Psi}
\Psi f(t)=\CL^{-1}_\beta[\Phi_\lambda \CL[f](\lambda)](t) \quad {\rm for} \quad f \in \CF^{-1}\CD(\R,X),
\end{equation}
where $\CL$ and $\CL^{-1}_\beta$ are the Laplace transform and its inverse defined in
\eqref{def:LT}.
Then 
\begin{equation} \label{est:maxreg}
\|e^{-\beta t}\Psi f\|_{L_p(\R,Y)} \leq C(p,X,Y)M\|e^{-\beta t}\|_{L_p(\R,X)} \quad \forall \; \beta \geq \beta_1.  
\end{equation}
\end{thm}
{\bf Proof:} For $\lambda=\beta+i\tau$ we have the following relation between 
Laplace and Fourier transforms defined in, respectively, \eqref{def:LT} and \eqref{def:FT}: 
\begin{align*}
&\CL[f](\lambda)=\int_\R e^{-\lambda t}f(t)dt = \CF[e^{-\beta t}f](\tau),\\
&\CL^{-1}_\beta[f](t)=\frac{1}{2\pi}\int_\R e^{\lambda t}f(\lambda)d\tau=e^{\beta t}\CF^{-1}[f](t).
\end{align*}
Therefore by \eqref{def:Psi} we have
$$
e^{-\beta t}\Psi f(t)=\CF^{-1}[\Psi_{\beta+i\tau}\CF[e^{-\beta t}](\tau)](t).
$$
Applying Theorem \ref{thm:Weis} to the above formula we conclude \eqref{est:maxreg}.

\qed

\noindent
{\bf Acknowledgements.} The work of Tomasz Piasecki was partially supported by National Science Centre grant
No2018/29/B/ST1/00339 (Opus). The work of Ewelina Zatorska was partially supported by the EPSRC Early Career Fellowship no. EP/V000586/1.

The authors would like to thank the anonymous Referees for numerous remarks which contributed significantly to the quality of the paper.

\section*{Declarations}

{\bf Conflict of interest:} the authors declare no conflict of interest.\\

\noindent
{\bf Data availability:} does not apply.


\begin{thebibliography}{99}
%
\bibitem{Ad} 
R. A. Adams, J. F. Fournier. 
\newblock \emph{Sobolev Spaces}, 
\newblock Second edition. Pure and Applied Mathematics (Amsterdam), 140. Elsevier/Academic Press, Amsterdam, (2003).

\bibitem{Bresch_chapter}
D. Bresch, B. Desjardins, J.-M. Ghidaglia, E. Grenier, M. Hilliairet. 
{\em Multifluid models including compressible fluids.} Handbook of Mathematical Analysis in Mechanics of Viscous Fluids, Giga, Y., Novotný, A. (Eds.) 1–52, (2017).

\bibitem{BDGG}
D. Bresch, B. Desjardins, J.-M. Ghidaglia, E. Grenier. 
{\em Global weak solutions to a generic two-fluid model.} Arch. Ration. Mech. Anal. 196, no. 2, 599–629, (2010).

\bibitem{BHL}
D.  Bresch, X. Huang, J. Li.
{\em Global weak solutions to one-dimensional non-conservative viscous compressible two-phase system.} 
Comm. Math. Phys. 309, no. 3, 737–755, (2012).

\bibitem{BMZ}
D. Bresch, P. B. Mucha, E. Zatorska. 
{\em Finite-energy solutions for compressible two-fluid Stokes system.} Arch. Ration. Mech. Anal. 232, no. 2, 987–1029,  (2019).

\bibitem{DHP} 
R.~Denk, M.~Hieber, J.~Pr\"u\ss.
\newblock {\em $\CR$-boundedness, Fourier multipliers and problems of elliptic and parabolic type}. 
\newblock Memoirs of AMS. Vol 166. no. 788,  (2003).

\bibitem{ES}
Y. Enomoto, Y. Shibata.
\newblock {\em On the $\CR$-sectoriality and the Initial Boundary Value Problem for the Viscous Compressible Fluid Flow}. 
\newblock Funkcialaj Ekvacioj, {\bf 56}, 441-505, (2013).

\bibitem{Evje1}
S. Evje, W. Wang, H. Wen. 
{\em Global well-posedness and decay rates of strong solutions to a non-conservative compressible two-fluid model.} 
Arch. Ration. Mech. Anal. 221, no. 3, 1285–1316, (2016).

\bibitem{Evje2}
S. Evje, H. Wen, C. Zhu. 
{\em On global solutions to the viscous liquid-gas model with unconstrained transition to single-phase flow.} 
Math. Models Methods Appl. Sci. 27, no. 2, 323–346, (2017).

\bibitem{Evje3}
Y. Qiao,  H. Wen, S. Evje. 
{\em Viscous two-phase flow in porous media driven by source terms: analysis and numerics.} SIAM J. Math. Anal. 51, no. 6, 5103–5140, (2019).

\bibitem{Guo}
 Z.H. Guo, J. Yang, L. Yao.
 {\em Global strong solution for a three-dimensional viscous liquid-gas two-phase flow model with vacuum.} J. Math. Phys., 52: 093102 (2011).

\bibitem{HM}
M. Hieber, M. Murata.
{\em The $L_p$-approach to the fluid-rigid body interaction problem for compressible fluids.} 
Evol. Equ. Control Theory 4 (2015), no. 1, 69--87

\bibitem{HNVW}
T. Hyt\"onen, J. van Neerven, M. Veraar, L. Weis. 
\newblock {\em Analysis in Banach spaces. Vol. I. Martingales and Littlewood-Paley theory.} A Series of Modern Surveys in Mathematics, vol. 63. Springer, Cham, (2016).

\bibitem{Ishii}
M. Ishii. 
{\em Thermo-Fluid Dynamic Theory of Two-Phase Flow.} 
Eyrolles, Paris, (1975).

\bibitem{JiNo}
B.J. Jin, A. Novotn\'y. 
{\it Weak-strong uniqueness for a bi-fluid model for a mixture of
non-interacting compressible fluids.} 
J. Differential Equations. 268, 204--238, (2019).

\bibitem{Kato}
T. Kato. 
\newblock {\em Perturbation Theory for Linear Operators.} 
\newblock Springer-Verlag Berlin Heidelberg, (1995).

\bibitem{KNP} O. Kreml, \v S. Ne\v casov\'a, T. Piasecki.
{\it Compressible Navier-Stokes system on a moving domain in the $L_p-L_q$ framework}. Waves in Flows, 127--158, in: Advances in Mathematical Fluid Mechanics, Birkh\"auser/Springer, 2021.


\bibitem{KSK} T. Kubo, Y. Shibata, K. Soga.
{\it On the $\CR$-boundedness for the two phase problem: compressible-incompressible model problem}
\newblock Boundary Value Problems, 141, (2014).

\bibitem{LiSuZa}
Y. Li, Y. Sun, E. Zatorska. 
{\it Large time bahavior for a compressible two-fluid model with algebraic pressure closure and large initial data.} Nonlinearity. 33, 4075-4094 (2020).

\bibitem{LiZa}
Y. Li, E. Zatorska. 
{\it On weak solutions to the compressible inviscid two-fluid model.} J. Differential Equations, 299, 33-50, (2021).


\bibitem{Murata} 
M.~Murata.
\newblock {\it On a maximal 
$L_p$-$L_q$ approach to the compressible viscous fluid flow
with slip boundary condition}.
\newblock Nonlinear Analysis, {106}, 86--109, (2014).


\bibitem{MS16}
 M.~Murata,  Y.~Shibata. 
 \newblock {\it On the global well-posedness for the compressible 
Navier-Stokes equations with slip boundary condition.}
\newblock J. Differential Equtions {\bf 260} (7), 5761--5795, (2016).


\bibitem{NoPo}
A. Novotn\'y, M. Pokorn\'y.
\newblock {\it Weak solutions for some compressible multicomponent 
fluid models.}
\newblock  Arch. Rational Mech. Anal. 235, 355--403, (2020).



\bibitem{PSZ} T. Piasecki, Y. Shibata, E. Zatorska.
\newblock {\em On strong dynamics of compressible two-component mixture flow}
\newblock  SIAM J. Math. Anal. 51, no. 4, 2793-2849, (2019).

\bibitem{PSZ2} T. Piasecki, Y. Shibata, E. Zatorska.
\newblock {\em On the isothermal compressible multi-component mixture flow: the local existence and maximal $L_p-L_q$ regularity of solutions.} 
\newblock Nonlinear Analysis 189, 111571, (2019).

\bibitem{PSZ3}
T. Piasecki, Y. Shibata, E. Zatorska
\newblock{\em On the maximal $L_p$-$L_q$  regularity of solutions
to a general linear parabolic system.}
\newblock  J. Differential Equations 268 (2020), no. 7, 3332--3369


\bibitem{S1} Y. Shibata.
\newblock {\em On the global well-posedness of some free boundary problem for a compressible barotropic viscous fluid flow.} 
\newblock Recent advances in partial differential equations and applications, 341-356, Contemp. Math., 666, Amer. Math. Soc., Providence, RI, (2016).

\bibitem{SS1} Y.~Shibata, S.~Shimizu.
\newblock {\it On some free boundary problem for the Navier-Stokes equtions.} 
\newblock Diff. Int. Eqns., {\bf 20}, 241--276, (2007).

\bibitem{SS2} 
Y.~Shibata, S. Shimizu.
{\it On the $L_p$-$L_q$ maximal regularity of the Neumann problem for  the Stokes equations in a bounded domain.} 
J. Reine Angew. Mat., {\bf 615}, 157--209, (2008).


\bibitem{Tanabe} 
H.~Tanabe. 
\newblock {\em Functional analytic methods for partial differential 
equations.} 
\newblock Monographs and textbooks in pure and
applied mathematics, Vol 204, Marchel Dekker, Inc. New York, Basel, 
(1997).

\bibitem{WaWeYa}

Y. Wang, H. Wen, L. Yao,
{\it On a Non-conservative Compressible Two-Fluid Model in a Bounded Domain: Global Existence and Uniqueness.}
 J. Math. Fluid Mech. 23, 4 (2021).

\bibitem{Weis} 
L.~Weis. 
\newblock {\em Operator-valued Fourier multiplier theorems and maximal $L_p$-regularity.} 
\newblock Math. Ann. 319, 735--758, (2001).

%
%
%
%
%
%
%
%
%
%
%
%
%
%
%
%
%
%
%
%
%
%
%
%
%
%
%
%
%
%
%
%
%
%
%
%
%
%
%
%
%
%
%
\end{thebibliography}
\end{document}